\documentclass[preprint,12pt]{elsarticle}
\usepackage{amsfonts}
\usepackage{amsmath}
\usepackage{mathrsfs}
\usepackage[applemac]{inputenc}
\usepackage{amssymb}
\usepackage{times, ulem, amsfonts}
\usepackage{mathpazo, amsmath, amssymb, amsthm,color}
\makeatletter
\def\ps@pprintTitle{%
   \let\@oddhead\@empty
   \let\@evenhead\@empty
   \let\@oddfoot\@empty
   \let\@evenfoot\@oddfoot
}
\makeatother

\def\q{{\mathcal Q}}

\def\p{{\mathcal P}}

\newcommand{\R}{{\mathbb{R}}}

\newcommand{\Z}{{\mathbb{Z}}}
\def\ddj{\dot \Delta_j}
\def\div{ \hbox{\rm div}\,  }

\newtheorem{theorem}{Theorem}[section]
\newtheorem{lemma}[theorem]{Lemma}
\newtheorem{definition}[theorem]{Definition}

\newtheorem{proposition}[theorem]{Proposition}
\newtheorem{remark}[theorem]{Remark}
\numberwithin{equation}{section}
\def\f{\frac}

\def\na{\nabla}

\def\ga{\Gamma}

\def\sa{\psi}
\def\aa{\varphi}

\def\La{\Lambda}

\let\wt=\widetilde

\begin{document}
\title{ {\Large\bf  Global  well-posedness to the  $n$-dimensional compressible Oldroyd-B model without damping mechanism}}

 \author{Xiaoping Zhai}
\ead{zhaixp@szu.edu.cn}

 \author{Zhi-Min Chen}
 \ead{zmchen@szu.edu.cn}
 \address{School  of Mathematics and Statistics, Shenzhen University, Shenzhen 518060, China}

\baselineskip=24pt

\begin{abstract}
The Cauchy problem of the compressible Oldroyd-B model without damping mechanism in $\mathbb{R}^n$ with $n\ge2$ is considered. The lack of  dissipation in density and stress tensor in the model
is compensated by exploiting an intrinsic structure  and introducing  new quantities between density, velocity and  stress tensor.
Therefore, global solutions  to the system with small initial data  in  critical  Besov spaces are obtained.
As a byproduct, optimal time decay rates  of the solutions are derived  by using an energy estimation argument.
The results remain  valid  for the compressible viscoelastic system  without the ``\!\!\div\!\!\!--\,{curl}'' structure assumption and thus  improve those given by Hu and   Wang  [ J. Differential Equations, {\bf 250},  1200--1231, 2011] and Qian and Zhang [Arch. Ration. Mech. Anal., {\bf 198}, 835--868, 2010].
\end{abstract}

\maketitle


\setcounter{theorem}{0}

\section{Introduction and the main results}
In the last several decades, non-Newtonian fluids which do not satisfy a linear relationship
between the stress tensor and the deformation tensor have been widely applied in engineering
and industry.
One particular subclass of non-Newtonian fluids is the  Oldroyd-B fluid,
which has been found to approximate the response to many dilute polymeric liquids.
Fluids of this type has a memory   and can
describe the motion of some viscoelastic flows, for example, the system coupling fluids and polymers.  Formulations about viscoelastic flows of Oldroyd-B type are first introduced by Oldroyd \cite{Oldroyd} and are extensively discussed in \cite{bird}.
In this paper, we 
 consider the wellposedness of the compressible Oldroyd-B model which has the following form \cite{Oldroyd}, \cite{tal}:
\begin{eqnarray}\label{sys1}
\left\{\begin{aligned}
&\partial_t\rho+\div(\rho u)=0,\\
&\rho(\partial_tu+(u\cdot\nabla) u)-\nu(\Delta u+\nabla\div u)+\nabla
P=\mu_1\div\tau,\\
&\partial_t\tau+(u\cdot\nabla)\tau+g(\tau,
\nabla u)+\beta\tau=\mu_2 D(u),\quad\quad\quad\quad x\in \R^n, \quad t>0.
\end{aligned}\right.
\end{eqnarray}
Here
 $\rho, u, \tau $ are the density,  velocity and
symmetric tensor of constrains, respectively.
$P=P(\rho)$, a smooth function of $\rho$,  is the pressure of  the fluid flow and  $g$ is a function defined as:
$$g(\tau, \nabla u)\stackrel{\mathrm{def}}{=}\tau
\Omega(u)-\Omega(u)\tau-b\left(D(u)\tau+\tau D(u)\right),\quad\hbox{with a parameter $b \in [-1,1] $}, $$
 The gradient $\nabla u$ is decomposed as the symmetric part and the skew-symmetric part expressed, respectively, as
\begin{equation}\label{omegau}
D(u) = \frac{1}{2} \big( \nabla u + (\nabla u)^{T} \big)\,\, \mbox{ and }\quad \Omega(u) = \frac{1}{2} \big( \nabla u - (\nabla u)^{T} \big).
\end{equation}
Moreover, we have the parameters   $\nu>0,\, \mu_1>0,\, \mu_2>0$ and $\beta\ge0.$
For  background details of  the
modeling, one may refer to   \cite{Bollada},  \cite{liuchun}, \cite{Oldroyd}, \cite{Schowalter} and references therein.

The theory of Oldroyd-B fluids attracts continuous attention of mathematicians. 
For the corresponding incompressible model ( $\rho=\mathrm{constant} $),
  the local and  global well-posedness in various
function spaces has  been studied by  Chen and Miao \cite{miaochangxing888},
Chemin and Masmoudi  \cite{chemin}, Elgindi and Liu  \cite{EL}, Elgindi and Rousset  \cite{ER},   Guillop\'{e} and Saut  \cite{GS, GS2}, Lin \cite{lin2012},  Renardy \cite{Renardy} and Zi {\it et al.}  \cite{zuiruizhao}.
 The global-in-time existence of weak solutions  was obtained by Lions and Masmoudi \cite{LM}
under the corotational derivative setting.
However, it's still open for the global-in-time existence of weak solutions in the general situation.

For the compressible model \eqref{sys1}, 
 Lei \cite{leizhen2006} and Gullop\'{e} {\it et al.} \cite {GST} investigated
its  incompressible limit problem 
 in  a
torus and a bounded domain of
$\R^3$, respectively.
Fang and Zi \cite{fang2014} further studied the incompressible limit problem in  $\R^n, n\ge2$, with ill prepared initial data in  Besov spaces.
Zi \cite{zuiruizhao2017} obtained the global small solutions of \eqref{sys1} in a critical $L^2$ Besov space.
Recently,
Barrett{\it et al.} \cite{Barrett} and Lu and Zhang  \cite{luyong} studied  the existence of global in-
time weak solutions in $\R^2$ and $\R^3$, respectively.
For the compressible Oldroyd type model
based on the deformation tensor, see the results  \cite{huxianpeng}, \cite{huxianpeng2013}, \cite{pan2019dcdsa}, \cite{pan2019},  \cite{qian} and references therein.

In the present paper, we consider   the global wellposedness of \eqref{sys1} without  damping mechanism  ($\beta=0$) in $\R^n (n\ge 2)$.
Without loss of generality, we  set the remaining  parameters  equal to  $1$ so  that  \eqref{sys1} is written as
\begin{eqnarray}\label{sys}
\left\{\begin{aligned}
&\partial_t\rho+\div(\rho u)=0,\\
&\rho(\partial_tu+(u\cdot\nabla) u)-(\Delta u+\nabla\div u)+\nabla
P=\div\tau,\\
&\partial_t\tau+(u\cdot\nabla)\tau+g(\tau,
\nabla u)= D(u),\quad\quad\quad\quad x\in \R^n, \quad t>0.
\end{aligned}\right.
\end{eqnarray}
The system \eqref{sys} is supplemented with the following initial conditions:
$$\rho|_{t=0}=\rho_0(x),\quad u|_{t=0}=u_0(x),\quad \tau|_{t=0}=\tau_0(x),\quad x\in \R^n,$$
and with far field behaviors
$$\rho\to 1,\quad\quad u\to0,\quad\quad \tau\to 0\quad \mathrm{as}
\quad|x|\to \infty.$$
Letting  $P'(1)=1,$
\begin{align}\label{}
 a\stackrel{\mathrm{def}}{=}  \rho-1, \quad     \aa\stackrel{\mathrm{def}}{=}P(1+a)-P(1), \quad k(a)\stackrel{\mathrm{def}}{=}(1+a){P'(1+a)}-1, \quad I(a)\stackrel{\mathrm{def}}{=}\frac{a}{1+a},
\end{align}
then we  rewrite \eqref{sys} into the following  system:
\begin{eqnarray}\label{m}
\left\{\begin{aligned}
&\partial_t \aa+u\cdot \nabla \aa + \div u =-k(a) \div u \, ,\\
&\partial_t\tau + u\cdot \nabla \tau   + g(\tau, \nabla u)  -  D (u)=0,\\
 &     \partial_t u + u\cdot\nabla u -\Delta u-\nabla\div u+ \nabla \aa -\div \tau
 = I(a)(\nabla \aa -\div \tau- \Delta u-\nabla\div u).
\end{aligned}\right.
\end{eqnarray}

To facilitate the understanding of the problem,
let us first recall  some known results for the
 Oldroyd-B model without  damping mechanism on the stress tensor.
 In fact, the most  of the  global solutions constructed in earlier studies  for the incompressible or compressible Oldroyd-B model depend heavily on the damping  term $\beta \tau$ in the third equation of \eqref{sys1}, since  the presence of the term ensures  $L^1$ or $L^2$ integration in time to be  available.  This  is crucial  to deal with the linear term $\div \tau$ in the momentum equation in  constructing  global small solutions.
For the incompressible  model without  the damping mechanism on the stress tensor, Zhu  \cite{zhuyi} obtained the global small  solutions in $\R^3,$ by
constructing two special time-weighted energies. This result was further extended to more general dimensions by  Chen and Hao  \cite{chenqionglei} and   Zhai \cite{zhaixiaopingdpde}
  in the critical $L^2$ Besov spaces and critical $L^p$ Besov spaces, respectively.
Recently, Zhu  \cite{zhuyi2} and Pan {\it et al.} \cite{pan2019} obtained the global small solutions to compressible
viscoelastic flows without the dvi-curl structure assumption  \cite[euqation (12)]{LeiLiuZhou} in $\R^3$  Sobolev spaces  and in $\R^n$  Besov  spaces, respectively.


However, it is still an open problem for the existence of   global solutions  to  the compressible Oldroyd-B model without  damping mechanism.
The main barrier lies in that we cannot get  the damping effect of the density and the stress tensor in the high frequencies part in the $L^p$ type Besov spaces. This  is different from the compressible Navier-Stokes equations treated in \cite{miaochangxing},  \cite{chenzhimin}, \cite{helingbing}, \cite{haspot}. The aim of the present paper is to break this barrier and solve this open problem by exploiting
an intrinsic structure of the model  and introducing several new quantities between the density, the velocity and the   stress tensor.

Let $\mathcal{S}(\R^n)$ be the space of
rapidly decreasing functions over $\R^n$ and $\mathcal{S}'(\R^n)$ be its dual
space. For any $z \in\mathcal{S}'(\R^n)$,
the lower and higher frequency parts are expressed as
\begin{align*}
z^\ell\stackrel{\mathrm{def}}{=}\sum_{j\leq j_0}\ddj z\quad\hbox{and}\quad
z^h\stackrel{\mathrm{def}}{=}\sum_{j>j_0}\ddj z
\end{align*}
for some fixed   integer $j_0\ge 1$ (the value of $j_0$ is dependent in  the proof of the main theorems). 
The corresponding  truncated semi-norms are defined  as follows:
\begin{align*}\|z\|^{\ell}_{\dot B^{s}_{p,1}}
\stackrel{\mathrm{def}}{=}  \|z^{\ell}\|_{\dot B^{s}_{p,1}}
\ \hbox{ and }\   \|z\|^{h}_{\dot B^{s}_{p,1}}
\stackrel{\mathrm{def}}{=}  \|z^{h}\|_{\dot B^{s}_{p,1}}.
\end{align*}

Denote $$
\Lambda\stackrel{\mathrm{def}}{=}\sqrt{-\Delta}, \quad\hbox{and}\quad \p=\mathcal{I}-\mathcal{Q}\stackrel{\mathrm{def}}{=}\mathcal{I}-\nabla\Delta^{-1}\div.$$

Now, we are in the position to  state the  first theorem of the paper:

\begin{theorem}\label{dingli}
Let   $n\ge 2$,  $(a_0^\ell,u_0,\tau_0^\ell)\in \dot{B}_{2,1}^{\frac n2-1}(\R^n)$ and  $(a^h_0,\tau_0^h)\in \dot{B}_{2,1}^{\frac n2}(\R^n)$.
 Assume the existence of  a positive constant $c_0$ such that,
\begin{align}\label{smallness}
\|(a^\ell_0,\tau_0^\ell)\|_{\dot{B}_{2,1}^{\frac {n}{2}-1}}+\| u_0\|_{\dot B^{\frac  n2-1}_{2,1}}+\|(a^h_0,\tau^h_0)\|_{\dot B^{\frac  n2}_{2,1}}\leq c_0.
\end{align}
Then
the system \eqref{m} has a unique global solution $(a,u,\tau)$ so that, for  any $T>0$,
\begin{align*}
&a^\ell\in C_b([0,T );{\dot{B}}_{2,1}^{\frac {n}{2}-1}(\R^n)),\quad a^h\in C_b([0,T );{\dot{B}}_{2,1}^{\frac {n}{2}}(\R^n)),
\\
& u\in C_b([0,T );{\dot{B}}_{2,1}^{\frac {n}{2}-1}\cap L^{1}
([0,T];{\dot{B}}_{2,1}^{\frac n2+1}(\R^n))\\
&\tau^\ell\in C_b([0,T );\dot{B}_{2,1}^{\frac n2-1}(\R^n)),\quad \tau^h\in C_b([0,T );\dot{B}_{2,1}^{\frac n2}(\R^n)), \\
&
(\Lambda^{-1}\p\div\tau)^\ell\in L^{1}
([0,T];{\dot{B}}_{2,1}^{\frac n2+1}(\R^n)),
 \quad (\Lambda^{-1}\p\div\tau)^h\in L^{1}
([0,T];{\dot{B}}_{2,1}^{\frac n2}(\R^n)).
\end{align*}
Moreover,  there exists some constant $C$ such that
\begin{align}\label{xiaonorm}
\|(a,\tau)\|&^\ell_{\widetilde{L}^\infty_t(\dot{B}_{2,1}^{\frac{n}{2}-1})}
+\|(a,\tau)\|^h_{\widetilde{L}^\infty_t(\dot{B}_{2,1}^{\frac{n}{2}})}+\|u\|_{\widetilde{L}^\infty_t(\dot{B}_{2,1}^{\frac{n}{2}-1})}+\|u\|_{L^1_t(\dot{B}_{2,1}^{\frac{n}{2}+1})}
\nonumber\\
&+\|(\Lambda^{-1}\p\div\tau)\|^\ell_{L^1_t(\dot{B}_{2,1}^{\frac{n}{2}+1})}
+\|\Lambda^{-1}\p\div\tau\|^h_{L^1_t(\dot{B}_{2,1}^{\frac{n}{2}})}\leq Cc_0.
\end{align}
\end{theorem}

\begin{remark}\label{23}
By using a  similar argument, one can  get the global small solution for the
compressible viscoelastic system  without any ``\!\!\div\!\!\!--\,{curl}'' structure assumption.
 Thus, Theorem \ref{dingli}  improves considerably the recent results by Pan and Xu  \cite{pan2019dcdsa} and Qian and Zhang   \cite{qian}.
\end{remark}

\begin{remark}\label{26}
Let $\rho$ be a constant, \eqref{sys1} reduces to the incompressible Oldroyd-B model without damping mechanism,
the above theorem coincides with the  Theorem 1.2 in \cite{chenqionglei} and  Theorem 1.2 ($p=2$) in \cite{zhaixiaopingdpde}.
\end{remark}

With the global solution shown above,
it is  natural to  explore  the large time asymptotic behavior of the solution.
The study of the large-time behavior of solutions is an important subject  in the field of  partial different equations. For example,   for   the compressible Navier-Stokes equations,
a celebrated  work in this subject    was  pioneered by  Matsumura and Nishida \cite{ma2} showing that  an initial function
 close to the steady-state solution $(\bar\rho,\bar u)=(1,0)$  in $H^3(\mathbb{R}^3)\times L^1(\mathbb{R}^3)$  develops a global solution behaving in    the asymptotic form
$$\|\nabla^k(\rho-\bar{\rho},u)(t)\|_{L^2}\le C(1+t)^{-\frac 34-\frac k2}\quad\mbox{for}\quad k=0,1.$$
Subsequently, Ponce \cite{ponce} obtained more general $L^p$ decay rates
$$\|\nabla^k(\rho-1,u)(t)\|_{L^p}\le C(1+t)^{-\frac n2(1-\frac 1p)-\frac k2},\quad2\leq p\leq\infty\quad 0\leq k\leq2 \quad n=2,3.$$
One can also refer to \cite{xujiang} for the decay rate of the compressible Navier-Stokes equations by  small perturbation in  Besov spaces.
Recently, Xin and Xu \cite{xujiang2019arxiv} obtained the decay rate of the compressible Navier-Stokes equations
without
the smallness of low frequencies of initial data.
However, there is almost no result  for the decay rate of the compressible Oldroyd-B model even with the damping mechanism.
Motivated by \cite{xujiang}, \cite{guoyan}  and  \cite{xujiang2019arxiv},  we give the  time-decay estimates of the global solutions constructed  in Theorem \ref{dingli}. Due to some technical reasons, here we only consider the pressure $P(\rho) $ being  a linear  function of $\rho$ in \eqref{m}. For the general  pressure function, some more complicated argument in Besov spaces is needed and is not  what we pursue in the present paper.

The  second theorem of the paper reads as follows:
\begin{theorem}\label{dingli2}
In addition to the assumptions of Theorem \ref{dingli},
let    $(a_0,u_0,\tau_0)\in{\dot{B}_{2,1}^{\sigma}}(\R^n)$  with
$-\frac n2<\sigma<\frac n2-1$.  Assume that the function $P(\rho)$ is linear. Then    the global small solution $(a,u,\tau)$  given by Theorem \ref{dingli} satisfies the decay property
\begin{eqnarray*}
\big\|\Lambda^{\alpha}(u,\Lambda^{-1}v)\big\|_{L^{q}} \le C(1+t)^{-\frac n4-\frac {(\alpha-\sigma)q-n}{2q}}
\end{eqnarray*}
for
$2\leq q\leq\infty$,
$\frac nq-\frac n2+\sigma<\alpha \leq\frac nq-1$ and   $v\stackrel{\mathrm{def}}{=}\nabla a-\div\tau$.
\end{theorem}

The rest of the paper unfolds as follows. In Section \ref{Section01},  we share the pivotal thought  of our analysis.
In Section \ref{Section2}, we shall collect some basic facts in Littlewood-Paley analysis and  product laws in
Besov spaces to be used for the proofs of the main theorems. In  Section \ref{Section3}, we will use four subsections to prove the main Theorem \ref{dingli}.
 Finally, in Section \ref{Section4}, we show  the proof  of Theorem \ref{dingli2}.

\section{Proving idea of the main theorems}\label{Section01}

The crucial part of the present paper is the proof of the first theorem.
We thus only  briefly describe the main difficulties arising from the derivation of  the global solution in Theorem \ref{dingli} and  the basic idea to overcome  them.

We start
 explaining the main ingredients of the proof of Theorem \ref{dingli}.
For simplicity, we only analysis the linearized equations of \eqref{sys}: 
\begin{eqnarray}\label{omas}
\left\{\begin{aligned}
&\partial_t a+ \div u =0 ,\\
&\partial_t\tau  - D (u)=0,\\
 &     \partial_t u  -(\Delta u+\nabla\div u)+ \nabla a -\div \tau=0.
\end{aligned}\right.
\end{eqnarray}
By neglecting  the effect of the stress tensor $\tau$, the system  \eqref{omas} reduces  to the  linearized compressible Navier-Stokes equations. It follows from   \cite{miaochangxing}, \cite{chenzhimin},  \cite{helingbing} and \cite{haspot}  that  the low frequency parts of the density and the  velocity have  smoothing effect in  $L^2$ type Besov spaces. Nevertheless, the high frequency part of the velocity field has  smoothing effect  and
 that of the density has  damping effect
in  $L^p$ type Besov spaces.

 When $a$ is a constant, equation \eqref{omas} becomes  the linearized incompressible Oldroyd-B model. From \cite{chenqionglei} and  \cite{zhaixiaopingdpde}, one can get  smoothing effect of $(\p u,\Lambda^{-1}\p\div\tau)$ in  low frequencies and smoothing effect of $\p u$ together with   damping effect of $\Lambda^{-1}\p\div\tau$ in  high frequencies.


For the compressible system   \eqref{omas}, one can follow the method for compressible Navier-Stokes equations to get  smoothing effect of $(a,\q u,\Lambda^{-1}\q\div\tau)$ in low frequencies. However, the  damping effect of $a$ and $\tau$  is not available in   high frequencies. Indeed, 
to overcome the difficulty, we rewrite the third equation of \eqref{omas} into the following form:
\begin{align*}
\partial_t \q u  -2\Delta   (\q u-\Delta^{-1} \nabla a +\Delta^{-1}\q\div \tau)=0.
\end{align*}
Just like the compressible Navier-Stoke equations, we
 introduce ``effective'' velocity field $$\Gamma=\q u-\Delta^{-1} \nabla a +\Delta^{-1}\q\div \tau,$$
    from which
    \begin{equation}\label{new1} \div \q u=\div \Gamma+  a-\Delta^{-1}\div(\q\div \tau),\quad \Lambda \q u=\Lambda \Gamma+ \Lambda ^{-1} \nabla a-\Lambda^{-1}(\q\div \tau).\end{equation}
Substituting  (\ref{new1}) into the first two equations of \eqref{omas}, we get
\begin{eqnarray}\label{omw}
\left\{\begin{aligned}
&\partial_t a+ a =-\div \Gamma- \Lambda^{-1}\div(\Lambda^{-1}\q\div \tau) ,\\
&\partial_t\Lambda^{-1}(\q\div \tau)  +\Lambda^{-1}(\q\div \tau)=-\Lambda \Gamma- \Lambda ^{-1} \nabla a.
\end{aligned}\right.
\end{eqnarray}
Due to the smoothing effect of $\Gamma$, the terms $\div \Gamma$ and $\Lambda \Gamma$ can be controlled easily. Yet, the
remaining  terms $ \Lambda^{-1}\div(\Lambda^{-1}\q\div \tau)$ and $  \Lambda ^{-1} \nabla a$ cannot be  controlled mutually,
since they have the same regularity in the high frequencies part. This difficult comes from  the occurrence of  three unknowns  in \eqref{omas}, which is different  to the  compressible Navier-Stokes equations involving two unknowns. That is, the damping effect of $a$ and $\tau$ in \eqref{omw} is not self-governed.
Thus, we cannot expect the availability of the damping effect of $a$ and $\tau$ in the high frequency part. This  gives rise to  a significant  difficulty to  bound the nonlinear term $\frac{1}{1+a}\nabla P(1+a)$ in the momentum equations
 and  thus some technique is required to get over it.
 In fact, from the above analysis of $a$ and $\tau$ in high frequencies, we find that the  new combination  $\Lambda^{-1}(\nabla a -\div \tau)$
  of $a$ and $\tau$ has  damping effect, although the effect  is absent   with respect to  the individual elements $a$ and  $\tau$.
Based  on this observation, we shall apply the new formulation  $\Lambda^{-1}(\nabla a -\div \tau)$  to close our energy estimates to be shown in the forthcoming   Section \ref{Section3}.

 In the proof of Theorem \ref{dingli}, 
 we shall employ the smoothing effect valid  for  the velocity field and the combination term  $\nabla a -\div \tau$. Hence, we  only expect to get the time decay of $(u,\nabla a -\div \tau)$  in the proof of Theorem \ref{dingli2}. The decay rate will be obtained in Section \ref{Section4} by solving a Lyapunov-type inequality which depends on the interpolation inequality in the low frequencies and high frequencies.


\setcounter{theorem}{0}
\setcounter{equation}{0}
\section{ Preliminaries }\label{Section2}

\noindent$\mathbf{Notation:}$ For operators  $A$ and  $B$, we denote  the commutator  $[A, B] = AB - BA$. The symbol  $a_1\lesssim a_2$ represents the inequality  $a_1 \le C a_2$ for a generic constant $C>0$. Given a Banach space $X$, we shall denote  $\|(f,g)\|_{X}=\|f\|_{X}+\|g\|_{X}$.

For  readers' convenience,  we list some basic knowledge on  Littlewood-Paley  theory.

\begin{definition}
 Let us consider smooth functions $\phi $ and $\chi$
on~$\R$  with the supports ${\rm supp} \phi \subset [\frac34,\frac83]$  and ${\rm supp} \chi \subset [0, 4/3]$ such that
$$
 \sum_{j\in\Z}\phi(2^{-j}s)=1 \,\,\mbox{ for }\,\, s>0
  \quad\hbox{and}\quad \chi(s)\stackrel{\mathrm{def}}{=} 1 - \sum_{j\geq
0}\phi(2^{-j}s) \,\,\mbox{ for }\,\, s \in \R. 
$$
Then we
define the blocks
$$
\ddj u={\mathcal F}^{-1}(\phi(2^{-j}|\xi|)\widehat{u})
 \quad\hbox{and}\quad \dot{S}_ju={\mathcal F}^{-1}(\chi(2^{-j}|\xi|)\widehat{u}).
$$
Let us remark that, for any homogeneous function $A$ of order 0  being smooth outside the origin  0. Then we have
\begin{equation*}\label{}
\|\ddj (A(D) f )\|_{L^p}\le C\|\ddj f\|_{L^p}\, \,\mbox{ for }\,\,  p\in[1,\infty].
\end{equation*}
Let $p, \, r \in [1,+\infty]$, $s\in \R$ and  $u\in\mathcal{S}'(\R^n)$. We define the Besov norm
$$
\|u\|_{{\dot{B}^s_{p,r}}}\stackrel{\mathrm{def}}{=}\big\|\big(2^{js}\|\ddj
u\|_{L^{p}}\big)_j\bigr\|_{\ell ^{r}({\mathop{\mathbb Z\kern 0pt}\nolimits})}
$$
and  then  the Besov space
$$\dot{B}_{p,r}^s(\R^n)\stackrel{\mathrm{def}}{=}\big\{u\in\mathcal{S}'_h(\R^n),\big|
\|u\|_{\dot{B}_{p,r}^s}<\infty\big\},$$
 where $\mathcal{S}'_h(\R^n)$ denotes the collection of  $u\in \mathcal{S}'(\R^n)$ so that  $\lim_{j\to -\infty}\|\dot{S}_ju\|_{L^\infty}=0$. 
\end{definition}

\begin{lemma}\label{bernstein}
Let $\mathcal{B}$ be a ball and $\mathcal{C}$ a ring of $\mathbb{R}^n$. For  $\lambda>0$,  integer $k\ge 0$,
$1\le p \le q\le\infty$ and a smooth homogeneous function $\sigma$ in  $\R^n\backslash\{0\}$ of degree m, then there hold
\begin{align*}
&& 
\|\nabla^{k}u\|_{L^q}\le C^{k+1}\lambda^{k+n(\frac1p-\frac1q)}\|u\|_{L^p}, \,\,\mbox{ whenever }\,\, \mathrm{supp} \,\hat{u}\subset\lambda \mathcal{B},\\
&& C^{-k-1}\lambda^k\|u\|_{L^p}\le \|\nabla^{k}u\|_{L^p}
\le C^{k+1}\lambda^{k}\|u\|_{L^p}, \,\,\mbox{ whenever }\,\,\mathrm{supp} \,\hat{u}\subset\lambda \mathcal{C},\\
&& \|\sigma(\nabla)u\|_{L^q}\le C_{\sigma,m}\lambda^{m+n(\frac1p-\frac1q)}\|u\|_{L^p}, \,\,\mbox{ whenever }\,\,\mathrm{supp} \,\hat{u}\subset\lambda \mathcal{C},
\end{align*}
where $\hat u$ denotes the Fourier transform of $u$.
\end{lemma}

Let us now state some classical
properties for the Besov spaces.
\begin{lemma}\label{qianru}

 Let  $1\le p\le \infty$, $0<\theta<1$, $s_2<s_1$ and  $u\in\dot{B}^{s_1}_{p,1}(\R^n)\cap\dot{B}^{s_2}_{p,1}(\R^n)$. Then there holds
\begin{align*}
\|u^\ell\|_{\dot{B}^{s_1}_{p,1}}\le C\| u^\ell\|_{\dot{B}^{s_2}_{p,1}}, \quad \|u^h\|_{\dot{B}^{s_2}_{p,1}}\le C\| u^h\|_{\dot{B}^{s_1}_{p,1}},
\quad \left[\dot{B}_{p,1}^{s_1}(\R^n),\dot{B}_{p,1}^{s_2}(\R^n)\right]_{\theta}=\dot{B}_{p,1}^{\theta
s_1+(1-\theta)s_2}(\R^n).
\end{align*}

Moreover, for a smooth homogeneous  function $\sigma$ in  $\R^n\backslash\{0\}$ of degree $m\in\Z$ , the operator $\sigma(\nabla)$ is a bounded operator mapping  $\dot{B}^{s}_{p,1}(\R^n)$ into $\dot{B}^{s-m}_{p,1}(\R^n)$.
\end{lemma}

In this paper, we frequently use the so-called ``time-space'' Besov spaces
or Chemin-Lerner space introduced by Chemin and Lerner \cite{bcd}.
\begin{definition}
Let $s\in\mathbb{R}$ and $0<T\leq +\infty$. We define
$$
\|u\|_{\widetilde{L}_{T}^{q}(\dot{B}_{p,1}^{s})}\stackrel{\mathrm{def}}{=}
\sum_{j\in\mathbb{Z}}
2^{js}\left(\int_{0}^{T}\|\dot{\Delta}_{j}u(t)\|_{L^p}^{q}dt\right)^{\frac{1}{q}}
$$
for  $ q,\,p \in [1,\infty)$ and with the standard modification for $p,\,q =\infty $.
\end{definition}
\noindent By Minkowski's inequality, we have the following inclusions between the
Chemin-Lerner space ${\widetilde{L}^\lambda_{T}(\dot{B}_{p,1}^s)}$ and the Bochner space ${{L}^\lambda_{T}(\dot{B}_{p,1}^s)}$:
\begin{align*}
\|u\|_{\widetilde{L}^\lambda_{T}(\dot{B}_{p,1}^s)}\le\|u\|_{L^\lambda_{T}(\dot{B}_{p,1}^s)}\hspace{0.5cm} \mathrm{if }\hspace{0.2cm}  \lambda\le r,\hspace{0.5cm}
\|u\|_{\widetilde{L}^\lambda_{T}(\dot{B}_{p,1}^s)}\ge\|u\|_{L^\lambda_{T}(\dot{B}_{p,1}^s)}\hspace{0.5cm} \mathrm{if }\hspace{0.2cm}  \lambda\> r.
\end{align*}

To study the pointwise function multiplications for the estimates of nonlinear terms of fluid motion equations, we need para-differential decomposition of Bony in the homogeneous context:
\begin{align}\label{bony}
uv=\dot{T}_uv+\dot{T}_vu+\dot{R}(u,v),
\end{align}
where
$$\dot{T}_uv\stackrel{\mathrm{def}}{=}\sum_{j\in Z}\dot{S}_{j-1}u\dot{\Delta}_jv, \hspace{0.5cm}\dot{R}(u,v)\stackrel{\mathrm{def}}{=}\sum_{j\in Z}
\dot{\Delta}_ju\widetilde{\dot{\Delta}}_jv,
\quad \widetilde{\dot{\Delta}}_jv\stackrel{\mathrm{def}}{=}\sum_{|j-j'|\le1}\dot{\Delta}_{j'}v.$$

The paraproduct $\dot{T}$ and the remainder $\dot{R}$ operators satisfy the following
continuous properties.

\begin{lemma}[\cite{bcd}]\label{fangji}
For  $s\in\mathbb{R}$, $\sigma\ge0$ and $1\leq p, \,p_1,\, p_2\leq\infty$  with $\frac{1}{p}=\frac{1}{p_1}+\frac{1}{p_2}$, then the
paraproduct $\dot T$ is a bilinear and  continuous operator mapping  $\dot{B}_{p_1,1}^{-\sigma}(\R^n)\times \dot{B}_{p_2,1}^s(\R^n)$ into
$\dot{B}_{p,1}^{s-\sigma}(\R^n)$  The remainder $\dot R$ is a bilinear continuous operator  mapping
$\dot{B}_{p_1, 1}^{s_1}(\R^n)\times \dot{B}_{p_2,1}^{s_2}(\R^n)$ to $
\dot{B}_{p,1}^{s_1+s_2}(\R^n)$ with
$s_1+s_2>0$.
\end{lemma}

In view of \eqref{bony}, Lemmas  \ref{bernstein} and  \ref{fangji}, one easily deduces the following  product laws:
\begin{lemma}\label{daishu}(see \cite[Proposition A.1]{xujiang})
Let $1\leq p, q\leq \infty$, $s_1\leq \frac{n}{q}$, $s_2\leq n\min\{\frac1p,\frac1q\}$ and $s_1+s_2>n\max\{0,\frac1p +\frac1q -1\}$. For $ (u,v)\in\dot{B}_{q,1}^{s_1}({\mathbb R} ^n)\times\dot{B}_{p,1}^{s_2}({\mathbb R} ^n)$, we have
\begin{align*}
\|uv\|_{\dot{B}_{p,1}^{s_1+s_2 -\frac{n}{q}}}\lesssim \|u\|_{\dot{B}_{q,1}^{s_1}}\|v\|_{\dot{B}_{p,1}^{s_2}}.
\end{align*}
\end{lemma}

The following  commutator estimate proved in
 \cite{bcd}
 will be employed  in our analysis.
\begin{lemma}(see   \cite[Lemma 2.100]{bcd})\label{jiaohuanzi}
Let $n\ge 2$, $1\leq p, q\leq \infty$,
$v \in \dot{B}_{q,1}^{s}(\R^n)$ and $\nabla u\in \dot{B}_{p,1}^{\frac {n}{p}}(\R^n)$.
Assume that
$$-n\min\left\{\frac 1p,1-\frac {1}{q}\right\}<s\leq 1+n\min\left\{\frac 1p,\frac 1q\right\}.$$
Then there holds the commutator estimate
$$
\|[\dot{\Delta}_j, u\cdot \nabla]v \|_{L^q}\lesssim d_j 2^{-js}\|\nabla u\|_{\dot{B}_{p,1}^{\frac {n}{p}}}\|v \|_{\dot{B}_{q,1}^{s}}.
$$
\end{lemma}

System \eqref{m} also involves compositions of functions $I(a)$  and $k(a)$, which are estimated by using  the following classical result:

  \begin{proposition}\label{fuhe} (\cite{bcd})
Let $F:\R\mapsto \R$ be a smooth function with $F(0)=0$,  $1\leq p$, $r\leq\infty$ and $s>0$. Then
$F: \dot {B}^{s}_{p,r}(\R^n)\cap L^{\infty}(\R^n)\mapsto  \dot {B}^{s}_{p,r}(\R^n)\cap L^{\infty}(\R^n)$  and
$$\left\|F(u)\right\|_{\dot B^{s}_{p,r}}\leq C\left\|u\right\|_{\dot B^{s}_{p,r}}$$
with  $C$ a constant depending only on $\left\|u\right\|_{L^{\infty}}$,   $s$, $p$,  $n$ and derivatives of $F$.

If $s>-\min\left(\frac {n}{p},\frac {n}{p'}\right)$,  then $F:\dot {B}^{s}_{p,r}(\R^n)\cap\dot {B}^{\frac {n}{p}}_{p,1}(\R^n)\mapsto  \dot {B}^{s}_{p,r}(\R^n)\cap\dot {B}^{\frac {n}{p}}_{p,1}(\R^n)$, and
\begin{equation*}
\left\|F(u)\right\|_{\dot B^{s}_{p,r}}\leq C(1+\|u\|_{\dot {B}^{\frac {n}{p}}_{p,1}})\|u\|_{\dot {B}^{s}_{p,r}}.
\end{equation*}

\end{proposition}

\setcounter{equation}{0}
\section{The proof of Theorem \ref{dingli}}\label{Section3}
To derive the small  global  solutions  to the Cauchy
problem \eqref{m},  it is convenient to  employ a continuity argument by extend a  local solution   to a global solution, if
 {\it a priori} solution estimates are ready. The analysis on the  local solution existence is standard. We thus  omit the process of deriving  the local solution.
In the present  paper,  we mainly   focus on the establishment of the {\it a priori} solution estimates.

 This section is divided  into four subsections to get the {\it a priori} estimates.  In the first  subsection, we get  the $L^\infty$ estimates of $(\aa,u,\tau)$ in  low frequencies. Then, we obtain the smoothing effect of $u$ and the hidden dissipation of $v=\nabla \aa -\div \tau$ in  low frequencies by introducing two new quantities in the second subsection.
In the third subsection, we follow  the approach in \cite{haspot} and introduce so-called ``effective'' velocity field
to find the  smoothing effect of $u$ and the damping effect of $v$ in  high frequencies.
 In the last subsection, the global solution is then obtained through a standard  continuity argument.

Note that $\aa=P(1+a)-P(1)$ with $P'(1)=1$. When $\|a\|_{L^\infty([0,T];\R^n)}\le\epsilon_0$ for $\epsilon_0$ small,
 the density difference  $a$ can be expressed as  a smooth function of $\aa$. That is, for convenience,   $a=\sa(\aa)$.
Throughout the proof, we assume  that
\begin{equation}\label{axiao}
\sup_{t\in\R_+,\, x\in\R^n} |a(t,x)|\leq \frac12,
\end{equation}
which  enables us to use freely the composition estimate stated in Proposition \ref{fuhe}.
The condition \eqref{axiao} is  ensured by the
small solution bound in \eqref{xiaonorm}, since  $\dot B^{\frac n2}_{2,1}(\R^n)\hookrightarrow L^\infty(\R^n).$

\subsection{The $L^\infty$ estimates of $(\aa,u,\tau)$ in  low frequencies}
Applying  the operator $\dot{\Delta}_{j}$ to (\ref{m}) and then using the standard energy argument, we have the following three equalities:
\begin{align*}
 &\f12\f{d}{dt} \|\ddj \aa\|_{L^2}^2 +\int\ddj \div u\cdot \ddj \aa{\,dx} \nonumber\\
 &\quad=-\int u\cdot \nabla\cdot \ddj \aa\cdot \ddj \aa{\,dx}- \int[\ddj ,u\cdot \na ]  \aa \cdot \ddj \aa{\,dx}-\int\ddj(k(a) \div u) \cdot \ddj \aa{\,dx},\nonumber\\
 &\f12\f{d}{dt} \|\ddj \tau\|_{L^2}^2 -\int\ddj D (u)\cdot \ddj \tau{\,dx} \nonumber\\
 &\quad=-\int u\cdot \nabla\cdot \ddj\tau \cdot \ddj \tau{\,dx}- \int[\ddj ,u\cdot \na ]  \tau \cdot \ddj \tau{\,dx}-\int\ddj(g(\tau, \nabla u)) \cdot \ddj \tau{\,dx},\nonumber\\
 &\f12\f{d}{dt} \|\ddj u\|_{L^2}^2 +\int|\ddj\nabla u|^2{\,dx} +\int|\ddj \div u|^2{\,dx}+\int\ddj\nabla \aa\cdot\ddj u \,dx-\int\ddj\div\tau \cdot\ddj u \,dx\nonumber\\
 &\quad=-\int u\cdot \nabla\cdot \ddj u \cdot \ddj u{\,dx}- \int[\ddj ,u\cdot \na ]  u \cdot \ddj u{\,dx}
\nonumber\\
 &\quad\quad+\int\ddj(I(a)({v}- \Delta u-\nabla\div u)) \cdot \ddj u{\,dx}.
 \end{align*}
 Integrating the previous equations from $0$ to $t$,   multiplying them by $2^{(\frac n2-1)j}$, and then  summing up the resultant equations for  $j\le j_0$, we get  that
\begin{align}\label{ert0}
&\|(\aa,u,\tau)^\ell\|_{\widetilde{L}_t^{\infty}(\dot{B}_{2,1}^{\frac n2-1})}
\lesssim\|(\aa_0,u_0,\tau_0)^\ell\|_{\dot{B}_{2,1}^{\frac {n}{2}-1}}
+ \int^t_0\|\nabla  u\|_{L^\infty}\|(\aa,u,\tau)^\ell\|_{\dot{B}_{2,1}^{\frac n2-1}}\,ds\nonumber\\
&\quad+ \int^t_0\|(k(a) \div u)^\ell\|_{\dot{B}_{2,1}^{\frac n2-1}}\,ds+ \int^t_0\|(I(a){v})^\ell\|_{\dot{B}_{2,1}^{\frac n2-1}}\,ds+ \int^t_0\|(g(\tau, \nabla u))^\ell\|_{\dot{B}_{2,1}^{\frac n2-1}}\,ds
\nonumber\\
&\quad+ \int^t_0\|(I(a)(\Delta u+\nabla\div u))^\ell\|_{\dot{B}_{2,1}^{\frac n2-1}}\,ds+\int^t_0\big(\sum_{j\le j_0}2^{(\frac n2-1)j}\| [\ddj ,u\cdot \na ] (\aa,u,\tau)\|_{L^2}\big),\,ds
\end{align}
where we have used by the H\"older inequality and the following cancellation identities:
\begin{align*}
&\int\ddj\div u\cdot\ddj \aa \,dx+\int\ddj\nabla \aa\cdot\ddj u \,dx=0,\quad\int\ddj D (u)\cdot\ddj \tau \,dx+\int\ddj\div\tau \cdot\ddj u \,dx=0.
\end{align*}

By virtue of embedding relation $\dot{B}_{2,1}^{\frac n2}(\R^n)\hookrightarrow L^\infty(\R^n),$ the integrand of the second term on the right-hand side of \eqref{ert0} can be bounded by
\begin{align}\label{yiming5+1+0}
\|\nabla  u\|_{L^\infty}\|(\aa,u,\tau)^\ell\|_{\dot{B}_{2,1}^{\frac n2-1}}
\lesssim&\|\nabla  u\|_{\dot{B}_{2,1}^{\frac n2}}\|(\aa,u,\tau)^\ell\|_{\dot{B}_{2,1}^{\frac n2-1}}
\lesssim\|(\aa,u,\tau)^\ell\|_{\dot{B}_{2,1}^{\frac n2-1}}\|u\|_{\dot{B}_{2,1}^{\frac n2+1}}.
\end{align}
By Lemma \ref{daishu}, the integrand of the third term on the right-hand side of \eqref{ert0} is estimated as
\begin{align}\label{yiming5+1-1-1}
\|(k(a) \div u)^\ell\|_{\dot{B}_{2,1}^{\frac n2-1}}
\lesssim& \|k(a)\|_{\dot{B}_{2,1}^{\frac n2-1}}\|\div u\|_{\dot{B}_{2,1}^{\frac n2}}
\lesssim \|k(a)\|_{\dot{B}_{2,1}^{\frac n2-1}}\|u\|_{\dot{B}_{2,1}^{\frac n2+1}}.
\end{align}
For bounding $\|k(a)\|_{\dot{B}_{2,1}^{\frac n2-1}},$  one cannot use the usual function compositions in Besov spaces directly, as $ n/2-1$ may be equal to zero in $\R^2$. Here, to lift this barrier, we  write
$$
k(a)=k'(0)\,a+a\wt k(a)\quad\hbox{with }\ \wt k(0)=0.
$$
Then  by Lemma \ref{daishu},  Proposition \ref{fuhe} and \eqref{axiao}, we  get
\begin{align}\label{wuwei}
 \|k(a)\|_{\dot{B}_{2,1}^{\frac n2-1}}
 \lesssim&k'(0)\|a\|_{\dot{B}_{2,1}^{\frac n2-1}}+\|a\wt k(a)\|_{\dot{B}_{2,1}^{\frac n2-1}}\nonumber\\
\lesssim&k'(0)\|a\|_{\dot{B}_{2,1}^{\frac n2-1}}+\|a\|_{\dot{B}_{2,1}^{\frac n2-1}}\|\wt k(a)\|_{\dot{B}_{2,1}^{\frac n2}}
 \nonumber\\
 \lesssim&k'(0)(\|a^\ell\|_{\dot{B}_{2,1}^{\frac n2-1}}+\|a^h\|_{\dot{B}_{2,1}^{\frac n2}})+(\|a^\ell\|_{\dot{B}_{2,1}^{\frac n2-1}}+\|a^h\|_{\dot{B}_{2,1}^{\frac n2}})(1+\|a\|_{L^\infty})\|a\|_{\dot{B}_{2,1}^{\frac n2}}
 \nonumber\\
 \lesssim&(\|a^\ell\|_{\dot{B}_{2,1}^{\frac n2-1}}+\|a^h\|_{\dot{B}_{2,1}^{\frac n2}}+1)(\|a^\ell\|_{\dot{B}_{2,1}^{\frac n2-1}}+\|a^h\|_{\dot{B}_{2,1}^{\frac n2}}).
\end{align}
From \eqref{wuwei}, equation \eqref{yiming5+1-1-1} is  further bounded by
\begin{align}\label{yiming5+1-1}
\|(k(a) \div u)^\ell\|_{\dot{B}_{2,1}^{\frac n2-1}}
\lesssim(\|\aa^\ell\|_{\dot{B}_{2,1}^{\frac n2-1}}+\|\aa^h\|_{\dot{B}_{2,1}^{\frac n2}}+1)(\|\aa^\ell\|_{\dot{B}_{2,1}^{\frac n2-1}}+\|\aa^h\|_{\dot{B}_{2,1}^{\frac n2}})\|u\|_{\dot{B}_{2,1}^{\frac n2+1}}.
\end{align}

The estimate \eqref{wuwei} remains  valid, if we  substitute $k(a)$ for $I(a)$. This helps to estimate the next integrand in \eqref{ert0}.
Thus repeating the same argument as above, we have
\begin{align}\label{yiming5+1+2}
\|(I(a) {v})^\ell\|_{\dot{B}_{2,1}^{\frac n2-1}}
\lesssim&\|I(a) {v}^\ell\|_{\dot{B}_{2,1}^{\frac n2-1}} +\|I(a) {v}^h\|_{\dot{B}_{2,1}^{\frac n2-1}} \nonumber\\
\lesssim&
\|I(a)\|_{\dot{B}_{2,1}^{\frac n2-1}}\|{v}^\ell\|_{\dot{B}_{2,1}^{\frac n2}}+\|I(a)\|_{\dot{B}_{2,1}^{\frac n2}}\|{v}^h\|_{\dot{B}_{2,1}^{\frac n2-1}}\nonumber\\
\lesssim& (\|\aa^\ell\|_{\dot{B}_{2,1}^{\frac n2-1}}+\|\aa^h\|_{\dot{B}_{2,1}^{\frac n2}}+1)(\|\aa^\ell\|_{\dot{B}_{2,1}^{\frac n2-1}}+\|\aa^h\|_{\dot{B}_{2,1}^{\frac n2}})\|{v}^\ell\|_{\dot{B}_{2,1}^{\frac n2}}\nonumber\\
&+(\|\aa^\ell\|_{\dot{B}_{2,1}^{\frac n2-1}}+\|\aa^h\|_{\dot{B}_{2,1}^{\frac n2}})\|{v}^h\|_{\dot{B}_{2,1}^{\frac n2-1}}.
\end{align}

For the estimate of  remaining terms in \eqref{ert0}, we  use  Lemma \ref{daishu} to obtain that
\begin{align}\label{yiming5+1+1}
\|(I(a)(\Delta u+\nabla\div u))^\ell\|_{\dot{B}_{2,1}^{\frac n2-1}}
\lesssim&\|I(a)\|_{\dot{B}_{2,1}^{\frac n2}}\|\nabla^2 u\|_{\dot{B}_{2,1}^{\frac n2-1}}
\lesssim(\|\aa^\ell\|_{\dot{B}_{2,1}^{\frac n2-1}}+\|\aa^h\|_{\dot{B}_{2,1}^{\frac n2}})\|u\|_{\dot{B}_{2,1}^{\frac n2+1}},\nonumber\\
\| (g(\tau, \nabla u))^\ell\|_{\dot{B}_{2,1}^{\frac n2-1}}
\lesssim&\|\tau\|_{\dot{B}_{2,1}^{\frac n2-1}}\|\nabla u\|_{\dot{B}_{2,1}^{\frac n2}}
\lesssim(\|\tau^\ell\|_{\dot{B}_{2,1}^{\frac n2-1}}+\|\tau^h\|_{\dot{B}_{2,1}^{\frac n2}})\|u\|_{\dot{B}_{2,1}^{\frac n2+1}}.
\end{align}
and, by Lemma \ref{jiaohuanzi}, 
\begin{align}\label{yiming5+1+5}
&\sum_{j\le j_0}2^{(\frac n2-1)j}\| [\ddj ,u\cdot \na ]   (\aa,u,\tau)\|_{L^2}\lesssim\big\|\nabla u\big\|_{\dot{B}^{\frac{n}{2}}_{2,1}}\big\|(\aa,u,\tau)\big\|_{\dot{B}^{\frac n2-1}_{2,1}}\nonumber\\
&\quad\lesssim(\big\|(\aa,u,\tau)^\ell\big\|_{\dot{B}^{ \frac n2-1}_{2,1}}+\big\|(\aa,u,\tau)^h\big\|_{\dot{B}^{\frac n2-1}_{2,1}})\big\|\nabla u\big\|_{\dot{B}^{\frac{n}{2}}_{2,1}}\nonumber\\
&\quad\lesssim(\big\|(\aa,u,\tau)^\ell\big\|_{\dot{B}^{ \frac n2-1}_{2,1}}+\big\|u^h\big\|_{\dot{B}^{\frac n2-1}_{2,1}}+\big\|(\aa,\tau)^h\big\|_{\dot{B}^{\frac n2}_{2,1}})(\|u^\ell\|_{\dot{B}_{2,1}^{\frac n2+1}}+\|u^h\|_{\dot{B}_{2,1}^{\frac n2+1}}).
\end{align}

Inserting \eqref{yiming5+1+0}, \eqref{yiming5+1-1}--\eqref{yiming5+1+5} into \eqref{ert0} gives the required estimate in low frequencies:
\begin{align}\label{yiming5+1+6}
&\|(\aa,u,\tau)^\ell\|_{\widetilde{L}_t^{\infty}(\dot{B}_{2,1}^{\frac n2-1})}
\lesssim\|(\aa_0,u_0,\tau_0)^\ell\|_{\dot{B}_{2,1}^{\frac {n}{2}-1}}+\int^t_0(\|\aa^\ell\|_{\dot{B}_{2,1}^{\frac n2-1}}+\|\aa^h\|_{\dot{B}_{2,1}^{\frac n2}})\|{v}^h\|_{\dot{B}_{2,1}^{\frac n2-1}}\,ds\nonumber\\
&\quad\quad+\int^t_0(\|\aa^\ell\|_{\dot{B}_{2,1}^{\frac n2-1}}+\|\aa^h\|_{\dot{B}_{2,1}^{\frac n2}}+1)(\|\aa^\ell\|_{\dot{B}_{2,1}^{\frac n2-1}}+\|\aa^h\|_{\dot{B}_{2,1}^{\frac n2}})(\|u\|_{\dot{B}_{2,1}^{\frac n2+1}}+\|{v}^\ell\|_{\dot{B}_{2,1}^{\frac n2}})\,ds\nonumber\\
&\quad\quad+ \int^t_0(\big\|(\aa,u,\tau)^\ell\big\|_{\dot{B}^{ \frac n2-1}_{2,1}}+\big\|u^h\big\|_{\dot{B}^{\frac n2-1}_{2,1}}+\big\|(\aa,\tau)^h\big\|_{\dot{B}^{\frac n2}_{2,1}})(\|u^\ell\|_{\dot{B}_{2,1}^{\frac n2+1}}+\|u^h\|_{\dot{B}_{2,1}^{\frac n2+1}})\,ds.
\end{align}

\subsection{The smoothing effect of $(u,\nabla \aa -\div \tau)$ in  low frequencies}
To  derive the smoothing effect of $(u,\nabla \aa -\div \tau)$ in  low frequencies,  we recall the new quantity
 ${v} 
 {=}\nabla \aa -\div \tau,$
and then  deduce from \eqref{m} that $({u},v)$ satisfies the following equations:
\begin{eqnarray}\label{ert2}
\left\{\begin{aligned}
&\partial_t{v} + u\cdot \nabla {v}   +\Delta u+\nabla\div u=f_1,\\
 &     \partial_t u + u\cdot\nabla u -(\Delta u+\nabla\div u)+ {v}= f_2,\\
\end{aligned}\right.
\end{eqnarray}
with
\begin{align*}
f_1\stackrel{\mathrm{def}}{=}\nabla u^{T}\nabla\tau-\nabla u^{T}\nabla \aa-\nabla(k(a) \div u)+\div g(\tau, \nabla u),\quad
f_2\stackrel{\mathrm{def}}{=}I(a)({v}- \Delta u-\nabla\div u).
\end{align*}

With  the aid of projection operators $\p$ and $\q$, we further separate the system \eqref{ert2} into the  compressible part
\begin{eqnarray}\label{ert4}
( \mathrm{Compressible} )\quad\left\{\begin{aligned}
&\partial_t\q{v}    +2\Delta \q u=\q f_1-\q(u\cdot\nabla{v}),\\
 &     \partial_t\q u  -2\Delta \q u+\q {v}=\q f_2 -\q(u\cdot\nabla u),\\
\end{aligned}\right.
\end{eqnarray}
and the incompressible part
\begin{eqnarray}\label{ert5}
( \mathrm{Incompressible} )\quad\left\{\begin{aligned}
&\partial_t\p{v} + u\cdot \nabla \p{v}   +\Delta \p u=\p f_1-[\p,u\cdot\nabla]{v},\\
 &     \partial_t\p u + u\cdot\nabla \p u -\Delta \p u+\p {v}= \p f_2-[\p,u\cdot\nabla]u.
\end{aligned}\right.
\end{eqnarray}
It's obvious that the systems  \eqref{ert4} and \eqref{ert5} have the same linear structure. Thus, we only pick out one of them, for example \eqref{ert5}, to reveal  the hidden dissipation of ${v}$.

For
$$\Lambda\stackrel{\mathrm{def}}{=}\sqrt{-\Delta},\quad  w\stackrel{\mathrm{def}}{=}\p u+\p{v},$$
we deduce from \eqref{ert5} that
\begin{eqnarray}\label{ert8}
\left\{\begin{aligned}
&\partial_t\La^{-1}\p{v}    -\La \p u=f_3,\\
&\partial_tw + u\cdot\nabla w+\p {v}=f_4
,\\
 &     \partial_t\p u + u\cdot\nabla \p u -\Delta \p u+\p {v}= f_5,
\end{aligned}\right.
\end{eqnarray}
with
\begin{align*}
f_3\stackrel{\mathrm{def}}{=}&\La^{-1}\p f_1-\La^{-1}\p(u\cdot\nabla{v}),\nonumber\\
f_4\stackrel{\mathrm{def}}{=}&\p f_1+\p f_2-[\p,u\cdot\nabla]{v}-[\p,u\cdot\nabla]u,\quad f_5\stackrel{\mathrm{def}}{=} \p f_2-[\p,u\cdot\nabla]u.
\end{align*}

Apply the operator $\dot{\Delta}_{j}$ to (\ref{ert8}). Then testing the resultant three  equations by $\ddj \La^{-1}\p{v}$,  $\ddj \p w$ and  $\ddj \p {u}$, respectively, we employ integration  by parts and  a standard commutator argument
to obtain
\begin{align}\label{ert11}
 &\f12\f{d}{dt} \|\ddj \La^{-1}\p{v}\|_{L^2}^2 -\int \ddj \p {v}\cdot \ddj \p u{\,dx} \lesssim \int |\ddj f_3  \cdot   \ddj \La^{-1}\p{v}  |{\,dx},
 \end{align}
\begin{align}\label{ert9}
 &\f12\f{d}{dt} \|\ddj w\|_{L^2}^2 +\|\ddj \p {v}  \|_{L^2}^2+\int\ddj \p {v}\cdot \ddj \p u{\,dx} \nonumber\\
 &\quad\lesssim \int |\div  \q u|  |\ddj w|^2{\,dx}+\int| [\ddj ,u\cdot \na ] w \cdot  \ddj w|{\,dx}   +\int |\ddj f_4  \cdot   \ddj w  |{\,dx},
 \end{align}
 and
\begin{align}\label{ert10}
 &\f12\f{d}{dt} \|\ddj \p u\|_{L^2}^2 +\|\ddj\nabla\p u  \|_{L^2}^2+\int \ddj \p {v}\cdot \ddj \p u{\,dx} \nonumber\\
 &\quad\lesssim \int |\div  \q u|  |\ddj \p u|^2{\,dx}+\int| [\ddj ,u\cdot \na ] \p u \cdot  \ddj \p u|{\,dx}   +\int |\ddj f_5  \cdot   \ddj \p u  |{\,dx}.
 \end{align}

Let $\eta \in(0,1)$ be a small  constant to be determined later. Summing up the equations $\eqref{ert9}\times\eta$, $\eqref{ert10}\times(1-\eta)$ and \eqref{ert11}, we have
\begin{align}\label{ert6}
 &\f12\f{d}{dt} (\|\ddj \La^{-1}\p{v}\|_{L^2}^2+\eta\|\ddj w\|_{L^2}^2 +(1-\eta)\|\ddj \p u\|_{L^2}^2)
 +(1-\eta)\|\ddj\nabla\p u  \|_{L^2}^2+\eta\|\ddj \p {v}  \|_{L^2}^2
 \nonumber\\
 &\quad\lesssim \int |\div  \q u|  (|\ddj w|^2,|\ddj \p u|^2){\,dx}+\int| [\ddj ,u\cdot \na ] (  w, \p u) \cdot  ( \ddj w,\ddj \p u)|{\,dx}\nonumber\\
 &\quad\quad   +\int |\ddj f_4  \cdot   \ddj w  |{\,dx}+\int |\ddj f_5  \cdot   \ddj \p u  |{\,dx}+\int |\ddj f_3  \cdot   \ddj \La^{-1}\p{v}  |{\,dx}.
 \end{align}
For the integer $j_0$ large enough and $j\le j_0$,  it holds that
\begin{align*}
\|  \ddj w\|_{L^2}\lesssim\|   \ddj \p u\|_{L^2}+\|   \ddj \p{v}\|_{L^2}
 \lesssim\|\ddj \p u\|_{L^2}+2^{j_0}\|\ddj \La^{-1}\p{v}\|_{L^2}.
\end{align*}
Thus, we  can find an $\eta>0$ small enough such that
\begin{eqnarray*}
&\|\ddj \La^{-1}\p{v}\|_{L^2}^2+\eta\|\ddj w\|_{L^2}^2 +(1-\eta)\|\ddj \p u\|_{L^2}^2
\geq \f1C( \|\ddj \p u\|_{L^2}^2+ \|\ddj \La^{-1}\p{v}\|_{L^2}^2).
\end{eqnarray*}%
Thanks to the H\"older inequality, integrating \eqref{ert6} over $[0, t]$,
 multiplying the resultant inequality by $2^{(\frac n2-1)j}$ and summing up about $j$ for $j\leq j_0$, we get that
\begin{align}\label{ert12}
&\|(\p u,\La^{-1}\p{v})^\ell\|_{\widetilde{L}_t^{\infty}(\dot{B}_{2,1}^{\frac n2-1})}
+\int^t_0\|(\p u,\La^{-1}\p{v})^\ell\|_{\dot{B}_{2,1}^{\frac n2+1}}\,ds\nonumber\\
&\quad\lesssim\|(\p u_0,\La^{-1}\p{v}_0)^\ell\|_{\dot{B}_{2,1}^{\frac {n}{2}-1}}+\int^t_0\sum_{j\le j_0}2^{(\frac n2-1)j}(\| [\ddj ,u\cdot \na ] ( \p u, \p{v})\|_{L^2})\,ds\nonumber\\
&\quad\quad
+ \int^t_0\|\nabla  u\|_{L^\infty}\|(\p u,\p{v})^\ell\|_{\dot{B}_{2,1}^{\frac n2-1}}\,ds+ \int^t_0\|(f_3,f_4,f_5)^\ell\|_{\dot{B}_{2,1}^{\frac n2-1}}\,ds.
\end{align}
This together with the inequality
 $\|([\p,u\cdot\nabla]{v})^\ell\|_{\dot{B}_{2,1}^{\frac n2-1}}\lesssim\|(\La^{-1}(u\cdot\nabla{v}))^\ell\|_{\dot{B}_{2,1}^{\frac n2-1}}
$
 implies that
\begin{align}\label{ert15-11}
&\|(\p u,\La^{-1}\p{v})^\ell\|_{\widetilde{L}_t^{\infty}(\dot{B}_{2,1}^{\frac n2-1})}
+\int^t_0\|(\p u,\La^{-1}\p{v})^\ell\|_{\dot{B}_{2,1}^{\frac n2+1}}\,ds\nonumber\\
&\quad\lesssim\|(\p u_0,\La^{-1}\p{v}_0)^\ell\|_{\dot{B}_{2,1}^{\frac {n}{2}-1}}+\int^t_0\sum_{j\le j_0}2^{(\frac n2-1)j}(\| [\ddj ,u\cdot \na ] ( \p u, \p{v})\|_{L^2})\,ds
\nonumber\\
&\quad\quad
+ \int^t_0\|\nabla  u\|_{L^\infty}\|(\p u,\p{v})^\ell\|_{\dot{B}_{2,1}^{\frac n2-1}}\,ds
+ \int^t_0\|( I(a)(\Delta u+\nabla\div u))^\ell\|_{\dot{B}_{2,1}^{\frac n2-1}}\,ds
\nonumber\\
&\quad\quad
+ \int^t_0\|(I(a){v})^\ell\|_{\dot{B}_{2,1}^{\frac n2-1}}\,ds+\int^t_0\|(k(a) \div u)^\ell\|_{\dot{B}_{2,1}^{\frac n2-1}}\,ds+\int^t_0\|(\div g(\tau, \nabla u))^\ell\|_{\dot{B}_{2,1}^{\frac n2-1}}\,ds\nonumber\\
&\quad\quad
+ \int^t_0\|([\p,u\cdot\nabla]u)^\ell\|_{\dot{B}_{2,1}^{\frac n2-1}}\,ds+ \int^t_0\|(\La^{-1}(u\cdot\nabla{v}))^\ell\|_{\dot{B}_{2,1}^{\frac n2-1}}\,ds\nonumber\\
&\quad\quad+\int^t_0\|(\La^{-1}(\nabla u^{T}\nabla \aa))^\ell\|_{\dot{B}_{2,1}^{\frac n2-1}}\,ds+\int^t_0\|(\La^{-1}(\nabla u^{T}\nabla\tau))^\ell\|_{\dot{B}_{2,1}^{\frac n2-1}}\,ds.
\end{align}

The nonlinear terms appeared on the right-hand  side of  \eqref{ert15-11} except the last four  are almost the same as those of \eqref{ert0}.
Hence, we are only concerned with these   four terms in the following.

Firstly, due to $$\ddj([\p,u\cdot\nabla]u)=[\ddj\p,u\cdot\nabla]u-[\ddj,u\cdot\nabla]\p u,$$ we can deduce from Lemma \ref{jiaohuanzi} that
\begin{align}\label{ert16}
\|([\p,u\cdot\nabla]u)^\ell\|_{\dot{B}_{2,1}^{\frac n2-1}}\lesssim(\big\|u^\ell\big\|_{\dot{B}^{ \frac n2-1}_{2,1}}+\big\|u^h\big\|_{\dot{B}^{\frac n2-1}_{2,1}})(\|u^\ell\|_{\dot{B}_{2,1}^{\frac n2+1}}+\|u^h\|_{\dot{B}_{2,1}^{\frac n2+1}}).
\end{align}

Secondly, we decompose $u\cdot{\nabla{v}}$ into $u\cdot{\nabla{v}^\ell}+u\cdot{\nabla{v}^h}$. Then from the product law in Besov spaces, we can readily  obtain, for any $n\ge 3$,
\begin{align}\label{ert16ewww}
\|(\La^{-1}({u}\cdot{\nabla{v}}))^\ell\|_{\dot{B}^{\frac n2-1}_{2,1}}
\lesssim&\|\La^{-1}({u}\cdot{\nabla{v}}^\ell)\|_{\dot{B}^{\frac n2-1}_{2,1}}+\|\La^{-1}({u}\cdot{\nabla{v}^h})\|_{\dot{B}^{\frac n2-1}_{2,1}}\nonumber\\
\lesssim&\|u\|_{\dot{B}_{2,1}^{\frac n2-1}}\|{v}^\ell\|_{\dot{B}_{2,1}^{\frac n2}}+\|u\|_{\dot{B}_{2,1}^{\frac n2}}\|{{v}}^h\|_{\dot{B}_{2,1}^{\frac n2-1}}.
\end{align}

However,
the above derivation is not valid for the  dimension $n=2$, since the function product
 operator  maps  $\dot{B}_{2,1}^{0}(\R^2)\times \dot{B}_{2,1}^{0}(\R^2)$ into $\dot{B}_{2,\infty}^{-1}(\R^2)$, but not into  $\dot{B}_{2,1}^{-1}(\R^2)$.
%
 To overcome this difficulty, we have to analysis carefully the frequency relations in the low part and the high part, respectively.

By using Bony's decomposition, one  write
\begin{align}\label{erw}
(\La^{-1}(u\cdot{\nabla{v}})=(\La^{-1}\dot{T}_{u}{\nabla{v}})+(\La^{-1}\dot{T}_{\nabla{v}} u)+(\La^{-1}\dot{R}(u,{\nabla{v}})).
\end{align}
Thanks to Bernstein's inequality and the H\"older inequality, the first term in \eqref{erw} is bounded as
\begin{align}\label{jiayou1}
&\|(\La^{-1}(\dot{T}_{u}{\nabla{v}}))^\ell\|_{\dot{B}^{\frac n2-1}_{2,1}}\nonumber\\
&\quad\lesssim \sum_{j\le j_0}2^{j(\frac n2-2)}\sum_{|j-k|\leq 4}\|\dot\Delta_j\big(\dot{S}_{k-1}u\dot{\Delta}_{k}{\nabla{v}}
\big)\|_{L^2}\nonumber\\
&\quad\lesssim \sum_{j\le j_0}2^{j(\frac n2-2)}\sum_{|j-k|\leq 4}  \Big(\sum_{k'\le k-2}\|\dot{\Delta}_{k'}u\|_{L^\infty}\Big)\|\dot{\Delta}_{k}{\nabla{v}}
\|_{L^2}\nonumber\\
&\quad\lesssim \sum_{j\le j_0}2^{j(\frac n2-2)}\Big(\sum_{\{|j-k|\leq 4\}\cap\{k\le j_0\}}  +\sum_{\{|j-k|\leq 4\}\cap\{k> j_0\}}  \Big)\Big(\sum_{k'\le k-2}\|\dot{\Delta}_{k'}u\|_{L^\infty}\Big)2^k\|\dot{\Delta}_{k}{\nabla{v}}
\|_{L^2}\nonumber\\
&\quad\lesssim\|u\|_{\dot{B}_{2,1}^{\frac n2-1}}\|{v}^\ell\|_{\dot{B}_{2,1}^{\frac n2}}+\|u\|_{\dot{B}_{2,1}^{\frac n2}}\|{v}^h\|_{\dot{B}_{2,1}^{\frac n2-1}}.
\end{align}

Similarly, we have
\begin{align}\label{jiayou2}
&\|(\La^{-1}(\dot{T}_{\nabla{v}}{u}))^\ell\|_{\dot{B}^{\frac n2-1}_{2,1}}
\lesssim \sum_{j\le j_0}2^{j(\frac n2-2)}\sum_{|j-k|\leq 4}\|\dot\Delta_j\big(\dot{S}_{k-1}{\nabla{v}}\dot{\Delta}_{k}u
\big)\|_{L^2}\nonumber\\
&\quad\lesssim \sum_{j\le j_0}2^{j(\frac n2-2)}\sum_{|j-k|\leq 4}  \Big(\sum_{k'\le k-2}\|\dot{\Delta}_{k'}{\nabla{v}}\|_{L^\infty}\Big)\|\dot{\Delta}_{k}u
\|_{L^2}\nonumber\\
&\quad\lesssim \sum_{j\le j_0}2^{j(\frac n2-2)}\sum_{|j-k|\leq 4} \Big(\big(\sum_{\{k'\le k-2\}\cap\{k'\le j_0\}}+\sum_{\{k'\le k-2\}\cap\{k'> j_0\}}\big)2^{(\frac n2+1)k'}\|\dot{\Delta}_{k'}{{v}}\|_{L^2}\Big)\|\dot{\Delta}_{k}u
\|_{L^2}\nonumber\\
&\quad\lesssim\|u\|_{\dot{B}_{2,1}^{\frac n2-1}}\|{v}^\ell\|_{\dot{B}_{2,1}^{\frac n2}}+\|u\|_{\dot{B}_{2,1}^{\frac n2}}\|{v}^h\|_{\dot{B}_{2,1}^{\frac n2-1}}.
\end{align}
Finally, to control the reminder term, we first recall $$\widetilde{\dot{\Delta}_k} \stackrel{\mathrm{def}}{=}\dot{\Delta}_k +\dot{\Delta}_{k-1}+\dot{\Delta}_{k+1}.$$
Thus  we have
\begin{align*}
\|(\La^{-1}\dot{R}(u,{\nabla{v}}))^\ell\|_{\dot{B}^{\frac n2-1}_{2,1}}\lesssim \mathcal{K}_1+\mathcal{K}_2+\mathcal{K}_3
\end{align*}
with
\begin{align*}
\mathcal{K}_1\stackrel{\mathrm{def}}{=}&\sum_{j\le j_0}2^{j(\frac n2-2)}\sum_{j-3\le k-1\le j}\| \dot{\Delta}_{k}u \|_{L^{ \infty}}\| \dot{\Delta}_{k-1} \nabla{v} \|_{L^{ 2}}\nonumber\\
&+\sum_{j\le j_0}2^{(n-1)j}\sum_{k-1>j}(\| \dot{\Delta}_{k}\La^{-1}u \|_{L^{ 2}}\|  \dot{\Delta}_{k-1}{\nabla{v}} \|_{L^{ 2}}+\| \dot{\Delta}_{k}u \|_{L^{2}}\| \dot{\Delta}_{k-1}{{v}} \|_{L^{ 2}})\nonumber\\
\stackrel{\mathrm{def}}{=}&\mathcal{K}_{1}^{(1)}+\mathcal{K}_{1}^{(2)},
\end{align*}
\begin{align*}
\mathcal{K}_2\stackrel{\mathrm{def}}{=}&\sum_{j\le j_0}2^{j(\frac n2-2)}\sum_{j-3\le k\le j}\| \dot{\Delta}_{k}u \|_{L^{ \infty}}\| \dot{\Delta}_{k} \nabla{v} \|_{L^{ 2}}\nonumber\\
&+\sum_{j\le j_0}2^{(n-1)j}\sum_{k>j}(\| \dot{\Delta}_{k}\La^{-1}u \|_{L^{ 2}}\|  \dot{\Delta}_{k}{\nabla{v}} \|_{L^{ 2}}+\| \dot{\Delta}_{k}u \|_{L^{2}}\| \dot{\Delta}_{k}{{v}} \|_{L^{ 2}})\nonumber\\
\stackrel{\mathrm{def}}{=}&\mathcal{K}_{2}^{(1)}+\mathcal{K}_{2}^{(2)},
\end{align*}
\begin{align*}
\mathcal{K}_3\stackrel{\mathrm{def}}{=}&\sum_{j\le j_0}2^{j(\frac n2-2)}\sum_{j-3\le k+1\le j}\| \dot{\Delta}_{k}u \|_{L^{ \infty}}\| \dot{\Delta}_{k+1} \nabla{v} \|_{L^{ 2}}\nonumber\\
&+\sum_{j\le j_0}2^{(n-1)j}\sum_{k+1>j}(\| \dot{\Delta}_{k}\La^{-1}u \|_{L^{ 2}}\|  \dot{\Delta}_{k+1}{\nabla{v}} \|_{L^{ 2}}+\| \dot{\Delta}_{k}u \|_{L^{2}}\| \dot{\Delta}_{k+1}{{v}} \|_{L^{ 2}})\nonumber\\
\stackrel{\mathrm{def}}{=}&\mathcal{K}_{3}^{(1)}+\mathcal{K}_{3}^{(2)}.
\end{align*}
By the Bernstein inequality, we infer that
\begin{align}\label{}
\mathcal{K}_{1}^{(1)}=&\sum_{j\le j_0}2^{j(\frac n2-2)}\sum_{j-3\le k-1\le j}\| \dot{\Delta}_{k}u \|_{L^{ \infty}}\| \dot{\Delta}_{k-1} \nabla{v} \|_{L^{ 2}}\nonumber\\
\lesssim&\sum_{j\le j_0}2^{j(\frac n2-2)}\sum_{j-3\le k-1\le j}2^{\frac{nk}{2}}\| \dot{\Delta}_{k}u \|_{L^{2}}2^k\| \dot{\Delta}_{k-1} {v} \|_{L^{ 2}}\nonumber\\
\lesssim&\|u\|_{\dot{B}_{2,1}^{\frac n2-1}}\|{v}^\ell\|_{\dot{B}_{2,1}^{\frac n2}}.
\end{align}
Similarly,
\begin{align}\label{}
\mathcal{K}_{2}^{(1)}+\mathcal{K}_{3}^{(1)}\lesssim\|u\|_{\dot{B}_{2,1}^{\frac n2-1}}\|{v}^\ell\|_{\dot{B}_{2,1}^{\frac n2}}.
\end{align}
The remaining terms $\mathcal{K}_{1}^{(2)}, \mathcal{K}_{2}^{(2)} $ and $\mathcal{K}_{3}^{(2)}$ have the same structure. We thus  only take a  representative  one as follows:
\begin{align}\label{}
\mathcal{K}_{1}^{(2)}=&\sum_{j\le j_0}2^{(n-1)j}\sum_{k-1>j}(\| \dot{\Delta}_{k}\La^{-1}u \|_{L^{ 2}}\|  \dot{\Delta}_{k-1}{\nabla{v}} \|_{L^{ 2}}+\| \dot{\Delta}_{k}u \|_{L^{2}}\| \dot{\Delta}_{k-1}{{v}} \|_{L^{ 2}})\nonumber\\
\lesssim&\|\La^{-1}u\|_{\dot{B}_{2,1}^{\frac n2+1}}\|\nabla{v}^h\|_{\dot{B}_{2,1}^{\frac n2-2}}+\|u\|_{\dot{B}_{2,1}^{\frac n2}}\|{v}^h\|_{\dot{B}_{2,1}^{\frac n2-1}}\lesssim2\|u\|_{\dot{B}_{2,1}^{\frac n2}}\|{v}^h\|_{\dot{B}_{2,1}^{\frac n2-1}}.
\end{align}
As a result, we get
\begin{align}\label{jiayou3}
&\|(\La^{-1}\dot{R}(u,{\nabla{v}}))^\ell\|_{\dot{B}^{\frac n2-1}_{2,1}}
\lesssim\|u\|_{\dot{B}_{2,1}^{\frac n2-1}}\|{v}^\ell\|_{\dot{B}_{2,1}^{\frac n2}}+\|u\|_{\dot{B}_{2,1}^{\frac n2}}\|{v}^h\|_{\dot{B}_{2,1}^{\frac n2-1}}.
\end{align}
The combination of  \eqref{jiayou1}, \eqref{jiayou2} and \eqref{jiayou3}  yields
\begin{align}\label{jiayou4}
\|(\La^{-1}({u}\cdot{\nabla{v}}))^\ell\|_{\dot{B}^{\frac n2-1}_{2,1}}\lesssim&\|u\|_{\dot{B}_{2,1}^{\frac n2-1}}\|{v}^\ell\|_{\dot{B}_{2,1}^{\frac n2}}+\|u\|_{\dot{B}_{2,1}^{\frac n2}}\|{v}^h\|_{\dot{B}_{2,1}^{\frac n2-1}}\nonumber\\
\lesssim&\|u\|_{\dot{B}_{2,1}^{\frac n2-1}}\|{v}^\ell\|_{\dot{B}_{2,1}^{\frac n2}}+\|u\|_{\dot{B}_{2,1}^{\frac n2-1}}^{\frac12}\|u\|_{\dot{B}_{2,1}^{\frac n2+1}}^{\frac12}\|{v}^h\|_{\dot{B}_{2,1}^{\frac n2-1}}^{\frac12}\|{v}^h\|_{\dot{B}_{2,1}^{\frac n2-1}}^{\frac12}\nonumber\\
\lesssim&\|u\|_{\dot{B}_{2,1}^{\frac n2-1}}\|{v}^\ell\|_{\dot{B}_{2,1}^{\frac n2}}+\|u\|_{\dot{B}_{2,1}^{\frac n2-1}}\|{v}^h\|_{\dot{B}_{2,1}^{\frac n2-1}}+\|{v}^h\|_{\dot{B}_{2,1}^{\frac n2-1}}\|u\|_{\dot{B}_{2,1}^{\frac n2+1}}.
\end{align}

Now, following the derivation of \eqref{jiayou4}, we have
\begin{align}\label{yiming3}
\|\La^{-1}(\nabla u^{T}\nabla \aa)^\ell\|_{\dot{B}_{2,1}^{\frac n2-1}}
\lesssim&(\|\aa^\ell\|_{\dot{B}_{2,1}^{\frac n2-1}}+\|\aa^h\|_{\dot{B}_{2,1}^{\frac n2}})(\|u^\ell\|_{\dot{B}_{2,1}^{\frac n2+1}}+\|u^h\|_{\dot{B}_{2,1}^{\frac n2+1}})
\end{align}
and
\begin{align}\label{yiming4}
\|\La^{-1}(\nabla u^{T}\nabla \tau)^\ell\|_{\dot{B}_{2,1}^{\frac n2-1}}
\lesssim&(\|\tau^\ell\|_{\dot{B}_{2,1}^{\frac n2-1}}+\|\tau^h\|_{\dot{B}_{2,1}^{\frac n2}})(\|u^\ell\|_{\dot{B}_{2,1}^{\frac n2+1}}+\|u^h\|_{\dot{B}_{2,1}^{\frac n2+1}}).
\end{align}

Inserting \eqref{yiming5+1+0}--\eqref{yiming5+1+5}, \eqref{ert16}, \eqref{ert16ewww}, \eqref{jiayou4}, \eqref{yiming3},  \eqref{yiming4} into \eqref{ert15-11} and combining the resultant equation with \eqref{yiming5+1+6}, we have
\begin{align}\label{yiming8}
&\|(\aa,u,\tau)^\ell\|_{\widetilde{L}_t^{\infty}(\dot{B}_{2,1}^{\frac n2-1})}
+\int^t_0\|(\p u,\La^{-1}\p{v})^\ell\|_{\dot{B}_{2,1}^{\frac n2+1}}\,ds\nonumber\\
&\quad\lesssim\|(\aa_0,u_0,\tau_0)^\ell\|_{\dot{B}_{2,1}^{\frac {n}{2}-1}}+\int^t_0(\|\aa^\ell\|_{\dot{B}_{2,1}^{\frac n2-1}}+\|\aa^h\|_{\dot{B}_{2,1}^{\frac n2}})\|{v}^h\|_{\dot{B}_{2,1}^{\frac n2-1}}\,ds\nonumber\\
&\quad\quad+\int^t_0(\|u\|_{\dot{B}_{2,1}^{\frac n2-1}}(\|{v}^\ell\|_{\dot{B}_{2,1}^{\frac n2}}+\|{v}^h\|_{\dot{B}_{2,1}^{\frac n2-1}})+\|{v}^h\|_{\dot{B}_{2,1}^{\frac n2-1}}\|u\|_{\dot{B}_{2,1}^{\frac n2+1}})\,ds\nonumber\\
&\quad\quad+\int^t_0(\|\aa^\ell\|_{\dot{B}_{2,1}^{\frac n2-1}}+\|\aa^h\|_{\dot{B}_{2,1}^{\frac n2}}+1)(\|\aa^\ell\|_{\dot{B}_{2,1}^{\frac n2-1}}+\|\aa^h\|_{\dot{B}_{2,1}^{\frac n2}})(\|u\|_{\dot{B}_{2,1}^{\frac n2+1}}+\|{v}^\ell\|_{\dot{B}_{2,1}^{\frac n2}})\,ds\nonumber\\
&\quad\quad+ \int^t_0(\big\|(\aa,u,\tau)^\ell\big\|_{\dot{B}^{ \frac n2-1}_{2,1}}+\big\|u^h\big\|_{\dot{B}^{\frac n2-1}_{2,1}}+\big\|(\aa,\tau)^h\big\|_{\dot{B}^{\frac n2}_{2,1}})(\|u^\ell\|_{\dot{B}_{2,1}^{\frac n2+1}}+\|u^h\|_{\dot{B}_{2,1}^{\frac n2+1}})\,ds.
\end{align}

As the linear structure of \eqref{ert4} is the same as \eqref{ert5},  the previous estimate \eqref{yiming8} with  $\|(\p u,\La^{-1}\p{v})^\ell\|_{\dot{B}_{2,1}^{\frac n2+1}}$
 replaced by $\|(\q u,\La^{-1}\q{v})^\ell\|_{\dot{B}_{2,1}^{\frac n2+1}}$ remains valid.

Thus, we  finally obtain the desired estimate
\begin{align}\label{yiming11}
&\|(\aa,u,\tau)^\ell\|_{\widetilde{L}_t^{\infty}(\dot{B}_{2,1}^{\frac n2-1})}
+\int^t_0\|(u,\La^{-1}{v})^\ell\|_{\dot{B}_{2,1}^{\frac n2+1}}\,ds\nonumber\\
&\quad\lesssim\|(\aa_0,u_0,\tau_0)^\ell\|_{\dot{B}_{2,1}^{\frac {n}{2}-1}}+\int^t_0(\|\aa^\ell\|_{\dot{B}_{2,1}^{\frac n2-1}}+\|\aa^h\|_{\dot{B}_{2,1}^{\frac n2}})\|{v}^h\|_{\dot{B}_{2,1}^{\frac n2-1}}\,ds\nonumber\\
&\quad\quad+\int^t_0(\|u\|_{\dot{B}_{2,1}^{\frac n2-1}}(\|{v}^\ell\|_{\dot{B}_{2,1}^{\frac n2}}+\|{v}^h\|_{\dot{B}_{2,1}^{\frac n2-1}})+\|{v}^h\|_{\dot{B}_{2,1}^{\frac n2-1}}\|u\|_{\dot{B}_{2,1}^{\frac n2+1}})\,ds\nonumber\\
&\quad\quad+\int^t_0(\|\aa^\ell\|_{\dot{B}_{2,1}^{\frac n2-1}}+\|\aa^h\|_{\dot{B}_{2,1}^{\frac n2}}+1)(\|\aa^\ell\|_{\dot{B}_{2,1}^{\frac n2-1}}+\|\aa^h\|_{\dot{B}_{2,1}^{\frac n2}})(\|u\|_{\dot{B}_{2,1}^{\frac n2+1}}+\|{v}^\ell\|_{\dot{B}_{2,1}^{\frac n2}})\,ds\nonumber\\
&\quad\quad+ \int^t_0(\big\|(\aa,u,\tau)^\ell\big\|_{\dot{B}^{ \frac n2-1}_{2,1}}+\big\|u^h\big\|_{\dot{B}^{\frac n2-1}_{2,1}}+\big\|(\aa,\tau)^h\big\|_{\dot{B}^{\frac n2}_{2,1}})(\|u^\ell\|_{\dot{B}_{2,1}^{\frac n2+1}}+\|u^h\|_{\dot{B}_{2,1}^{\frac n2+1}})\,ds.
\end{align}

\subsection{The estimates of $(\aa,u,\tau)$ in  high frequencies}
To  get the smoothing effect of $u$ and the damping effect of ${v}$ in  high frequencies, we will use  the  effective velocity field given in \cite{haspot}.

 In  the  Subsection 4.2 for the analysis in low frequencies, we have obtained  the smoothing effect of $v$ in the compressible part and  the incompressible part, respectively.
  In contrast, in the present section,  we take  the compressible part for example to get the smoothing effect of $u$ and the damping effect of ${v}$ in the high frequencies.

Letting
$$\ga_1\stackrel{\mathrm{def}}{=}\q u-\f12\Delta^{-1}\q {v},\quad \ga_2\stackrel{\mathrm{def}}{=}\p u-\Delta^{-1}\p {v},$$
we  deduce from \eqref{ert4} that
\begin{align}\label{yiming14}
\partial_t \ga_1-2\Delta \ga_1
=&\q u+\q f_2-\q(u\cdot\nabla u)-\Delta^{-1}(\q f_1-\q(u\cdot \nabla {v}) )\nonumber\\
=&\ga_1+\f12\Delta^{-1}\q {v}+\q f_2-\q(u\cdot\nabla u)-\Delta^{-1}(\q f_1-\q(u\cdot \nabla {v}) ).
\end{align}
This yields, by a standard energy argument,
\begin{align}\label{yiming15}
&\|\ga_1^h\|_{ \widetilde{L}_t^\infty(\dot B^{\frac  n2-1}_{2,1})}+ \|\ga_1^h\|_{L^1_t(\dot B^{\frac  n2+1}_{2,1})}\nonumber\\
&\quad\lesssim \|(\ga_1)_0^h\|_{\dot B^{\frac  n2-1}_{2,1}}+\int_0^t(\|\ga_1^h\|_{\dot B^{\frac  n2-1}_{2,1}}+\f12\|(\q {v})^h\|_{\dot B^{\frac  n2-3}_{2,1}})\,ds+\int_0^t\|(\q f_2,\Delta^{-1}\q f_1)^h\|_{\dot B^{\frac  n2-1}_{2,1}}\,ds\nonumber\\
&\quad\quad+\int_0^t(\|(\q(u\cdot\nabla u))^h\|_{\dot B^{\frac  n2-1}_{2,1}}+\|(\Delta^{-1}\q(u\cdot \nabla {v}))^h\|_{\dot B^{\frac  n2-1}_{2,1}})\,ds.
\end{align}

Plugging $\q u=\ga_1+\f12\Delta^{-1}\q {v}$ into the first equation in \eqref{ert4} gives
\begin{align}\label{yiming16}
&\partial_t \q{v}+u\cdot \nabla \q{v} +\q{v}=-2\Delta\ga_1+\q f_1-[\q,u\cdot \nabla]{v}.
\end{align}
       Applying the operator $\ddj$ onto  \eqref{yiming16},  then taking the inner product with  $\ddj \q{v}$ and  summing up  the resultant equations with respect to  the high frequency indices $j>j_0$, we have
\begin{align}\label{yiming17}
&\|(\q{v})^h\|_{ \widetilde{L}_t^\infty(\dot B^{\frac  n2-1}_{2,1})}+ \|(\q{v})^h\|_{L^1_t(\dot B^{\frac  n2-1}_{2,1})}\nonumber\\
&\quad\lesssim \|(\q{v})_0^h\|_{\dot B^{\frac  n2-1}_{2,1}}+\int_0^t(\|\ga_1^h\|_{\dot B^{\frac  n2+1}_{2,1}}+\|(\q f_1)^h\|_{\dot B^{\frac  n2-1}_{2,1}}+\|([\q,u\cdot \nabla]{v})^h\|_{\dot B^{\frac  n2-1}_{2,1}})\,ds
\nonumber\\
&\quad\quad+\int_0^t\|\div u\|_{L^\infty}\|(\q{v})^h\|_{\dot B^{\frac  n2-1}_{2,1}}\,ds+\int_0^t\sum_{j\ge j_0}2^{(\frac {n}{2}-1)j}\|[\ddj,u\cdot\nabla]\q{v}\|_{L^2}\,ds
.
\end{align}
Owing to the high frequency cut-off at $|\xi|\sim 2^{j_0},$ we have
$$
\|\ga_1\|^h_{L^1(\dot B^{\frac  n2-1}_{2,1})}\lesssim 2^{-2j_0}\|\ga_1\|^h_{L^1(\dot B^{\frac  n2+1}_{2,1})}
\quad\hbox{and}\quad
\|\q{v}\|^h_{L^1(\dot B^{\frac  n2-3}_{2,1})}\lesssim  2^{-2j_0}\|\q{v}\|^h_{L^1(\dot B^{\frac  n2-1}_{2,1})}.
$$
Multiplying   \eqref{yiming15} by a suitable large constant and then inserting the resultant equation into \eqref{yiming17}, we get
\begin{align}\label{yiming18}
&\|\ga_1^h\|_{ \widetilde{L}_t^\infty(\dot B^{\frac  n2-1}_{2,1})}+\|(\q{v})^h\|_{ \widetilde{L}_t^\infty(\dot B^{\frac  n2-1}_{2,1})}+ \|\ga_1^h\|_{L^1_t(\dot B^{\frac  n2+1}_{2,1})}+ \|(\q{v})^h\|_{L^1_t(\dot B^{\frac  n2-1}_{2,1})}\nonumber\\
&\quad\lesssim \|(\ga_1)_0^h\|_{\dot B^{\frac  n2-1}_{2,1}}+\|(\q{v})_0^h\|_{\dot B^{\frac  n2-1}_{2,1}}+\int_0^t(\|(\q f_1,\q f_2)^h\|_{\dot B^{\frac  n2-1}_{2,1}}+\|\div u\|_{L^\infty}\|(\q{v})^h\|_{\dot B^{\frac  n2-1}_{2,1}})\,ds\nonumber\\
&\quad\quad+\int_0^t\|([\q,u\cdot \nabla]{v})^h\|_{\dot B^{\frac  n2-1}_{2,1}}ds+\int_0^t\sum_{j\ge j_0}2^{(\frac {n}{2}-1)j}\|[\ddj,u\cdot\nabla]\q{v}\|_{L^2}\,ds\nonumber\\
&\quad\quad+\int_0^t(\|(\q(u\cdot\nabla u))^h\|_{\dot B^{\frac  n2-1}_{2,1}}+\|(\Delta^{-1}\q(u\cdot \nabla {v})^h\|_{\dot B^{\frac  n2-1}_{2,1}})\,ds.
\end{align}
By Lemma \ref{daishu}, we have
\begin{align}\label{shiping4}
\|(\q f_1)^h\|_{\dot B^{\frac  n2-1}_{2,1}}\lesssim&\|\nabla u^{T}\nabla \aa\|_{\dot B^{\frac  n2-1}_{2,1}}+\|\nabla u^{T}\nabla \tau\|_{\dot B^{\frac  n2-1}_{2,1}}\nonumber\\
&+\|\nabla(k(a) \div u)\|_{\dot B^{\frac  n2-1}_{2,1}}+\|\div g(\tau, \nabla u)\|_{\dot B^{\frac  n2-1}_{2,1}}\nonumber\\
\lesssim&(\|\nabla \aa\|_{\dot B^{\frac  n2-1}_{2,1}}+\|k(a)\|_{\dot B^{\frac  n2}_{2,1}}+\|\nabla \tau\|_{\dot B^{\frac  n2-1}_{2,1}}+\|\tau\|_{\dot B^{\frac  n2}_{2,1}})\|\nabla u\|_{\dot B^{\frac  n2}_{2,1}}\nonumber\\
\lesssim&(\|(\aa,\tau)^\ell\|_{\dot{B}_{2,1}^{\frac n2-1}}+\|(\aa,\tau)^h\|_{\dot{B}_{2,1}^{\frac n2}})(\|u^\ell\|_{\dot{B}_{2,1}^{\frac n2+1}}+\|u^h\|_{\dot{B}_{2,1}^{\frac n2+1}}).
\end{align}
Following the derivation of   \eqref{yiming5+1+2} and  \eqref{yiming5+1+1}, we obtain the estimate of $\q f_2$ as
\begin{align}\label{shiping66}
\|(\q f_2)^h\|_{\dot B^{\frac  n2-1}_{2,1}}\lesssim& (\|\aa^\ell\|_{\dot{B}_{2,1}^{\frac n2-1}}+\|\aa^h\|_{\dot{B}_{2,1}^{\frac n2}}+1)(\|\aa^\ell\|_{\dot{B}_{2,1}^{\frac n2-1}}+\|\aa^h\|_{\dot{B}_{2,1}^{\frac n2}})\|{v}^\ell\|_{\dot{B}_{2,1}^{\frac n2}}\nonumber\\
&+(\|\aa^\ell\|_{\dot{B}_{2,1}^{\frac n2-1}}+\|\aa^h\|_{\dot{B}_{2,1}^{\frac n2}})(\|{v}^h\|_{\dot{B}_{2,1}^{\frac n2-1}}+\|u\|_{\dot{B}_{2,1}^{\frac n2+1}}).
\end{align}
By the embedding relation $\dot{B}_{2,1}^{\frac n2}(\R^n)\hookrightarrow L^\infty(\R^n)$ and the inequality  $$\|{v}^h\|_{\dot{B}_{2,1}^{\frac n2-1}}\lesssim\|(\aa,\tau)^h\|_{\dot{B}_{2,1}^{\frac n2}},$$
we  get
\begin{align}\label{shiping8}
\|\div u\|_{L^\infty}\|(\q{v})^h\|_{\dot B^{\frac  n2-1}_{2,1}}\lesssim&\|{v}^h\|_{\dot B^{\frac  n2-1}_{2,1}}\|\nabla u\|_{\dot B^{\frac  n2}_{2,1}}\nonumber\\
\lesssim&\|{v}^h\|_{\dot{B}_{2,1}^{\frac n2-1}}(\|u^\ell\|_{\dot{B}_{2,1}^{\frac n2+1}}+\|u^h\|_{\dot{B}_{2,1}^{\frac n2+1}})\nonumber\\
\lesssim&\|(\aa,\tau)^h\|_{\dot{B}_{2,1}^{\frac n2}}(\|u^\ell\|_{\dot{B}_{2,1}^{\frac n2+1}}+\|u^h\|_{\dot{B}_{2,1}^{\frac n2+1}}).
\end{align}
Similarly, we have
\begin{align}\label{shiping11}
\|(\q(u\cdot\nabla u))^h\|_{\dot B^{\frac  n2-1}_{2,1}}
\lesssim&(\|u^\ell\|_{\dot{B}_{2,1}^{\frac n2-1}}+\|u^h\|_{\dot{B}_{2,1}^{\frac n2-1}})(\|u^\ell\|_{\dot{B}_{2,1}^{\frac n2+1}}+\|u^h\|_{\dot{B}_{2,1}^{\frac n2+1}}).
\end{align}
From the identity
$$
\ddj([\q,u\cdot \nabla]{v})=[\ddj\q,u\cdot \nabla]{v}-[\ddj,u\cdot \nabla]\q{v}
$$ and Lemma \ref{jiaohuanzi}, we have
\begin{align}\label{shiping10}
&\sum_{j\ge j_0}2^{(\frac {n}{2}-1)j}\|[\ddj,u\cdot\nabla]\q{v}\|_{L^2}+\|([\q,u\cdot \nabla]{v})^h\|_{\dot B^{\frac  n2-1}_{2,1}}\nonumber\\
&\quad\lesssim\|{v}\|_{\dot B^{\frac  n2-1}_{2,1}}\|\nabla u\|_{\dot B^{\frac  n2}_{2,1}}\nonumber\\
&\quad\lesssim(\|{v}^\ell\|_{\dot{B}_{2,1}^{\frac n2-2}}+\|{v}^h\|_{\dot{B}_{2,1}^{\frac n2-1}})(\|u^\ell\|_{\dot{B}_{2,1}^{\frac n2+1}}+\|u^h\|_{\dot{B}_{2,1}^{\frac n2+1}})\nonumber\\
&\quad\lesssim(\|(\aa,\tau)^\ell\|_{\dot{B}_{2,1}^{\frac n2-1}}+\|(\aa,\tau)^h\|_{\dot{B}_{2,1}^{\frac n2}})(\|u^\ell\|_{\dot{B}_{2,1}^{\frac n2+1}}+\|u^h\|_{\dot{B}_{2,1}^{\frac n2+1}}).
\end{align}
Using the equation $u\cdot \nabla {v}=\div(u\otimes{v})-{v}\div u$,  we  get
\begin{align}\label{shiping12}
&\|((\Delta^{-1}\q(u\cdot \nabla {v}))^h\|_{\dot B^{\frac  n2-1}_{2,1}}
\lesssim\|u {v}\|_{\dot B^{\frac  n2-1}_{2,1}}+\|{v}\div u \|_{\dot B^{\frac  n2-1}_{2,1}}\nonumber\\
&\quad\lesssim\|u \|_{\dot B^{\frac  n2}_{2,1}}\| {v}\|_{\dot B^{\frac  n2-1}_{2,1}}+\|{v} \|_{\dot B^{\frac  n2-1}_{2,1}}\|\div u \|_{\dot B^{\frac  n2}_{2,1}}\lesssim\|u \|_{\dot B^{\frac  n2}_{2,1}}^2+\| {v}\|_{\dot B^{\frac  n2-1}_{2,1}}^2+\|{v} \|_{\dot B^{\frac  n2-1}_{2,1}}\| u\|_{\dot B^{\frac  n2+1}_{2,1}}\nonumber\\
&\quad\lesssim\|u^\ell \|_{\dot B^{\frac  n2-1}_{2,1}}\|u^\ell \|_{\dot B^{\frac  n2+1}_{2,1}}+\|u^h \|_{\dot B^{\frac  n2-1}_{2,1}}
\|u^h \|_{\dot B^{\frac  n2+1}_{2,1}}+\| {v}^\ell\|_{\dot B^{\frac  n2-2}_{2,1}}\| {v}^\ell\|_{\dot B^{\frac  n2}_{2,1}}\nonumber\\
&\quad\quad+\| {v}^h\|_{\dot B^{\frac  n2-1}_{2,1}}\| {v}^h\|_{\dot B^{\frac  n2-1}_{2,1}}
+(\|{v}^\ell \|_{\dot B^{\frac  n2-2}_{2,1}}+\|{v}^h \|_{\dot B^{\frac  n2-1}_{2,1}})(\| u^\ell \|_{\dot B^{\frac  n2+1}_{2,1}}+\|u^h \|_{\dot B^{\frac  n2+1}_{2,1}})\nonumber\\
&\quad\lesssim\|u^\ell \|_{\dot B^{\frac  n2-1}_{2,1}}\|u^\ell \|_{\dot B^{\frac  n2+1}_{2,1}}+\|u^h \|_{\dot B^{\frac  n2-1}_{2,1}}
\|u^h \|_{\dot B^{\frac  n2+1}_{2,1}}\nonumber\\
&\quad\quad+
(\|(\aa,\tau)^\ell\|_{\dot{B}_{2,1}^{\frac n2-1}}+\|(\aa,\tau)^h\|_{\dot{B}_{2,1}^{\frac n2}})(\| u^\ell \|_{\dot B^{\frac  n2+1}_{2,1}}+\| {v}^\ell\|_{\dot B^{\frac  n2}_{2,1}}+\|u^h \|_{\dot B^{\frac  n2+1}_{2,1}}+\| {v}^h\|_{\dot B^{\frac  n2-1}_{2,1}})
.
\end{align}

Inserting \eqref{shiping4}, \eqref{shiping66}--\eqref{shiping12} into \eqref{yiming18} gives
\begin{align}\label{shiping13}
&\|\ga_1^h\|_{ \widetilde{L}_t^\infty(\dot B^{\frac  n2-1}_{2,1})}+\|(\q{v})^h\|_{ \widetilde{L}_t^\infty(\dot B^{\frac  n2-1}_{2,1})}+ \|\ga_1^h\|_{L^1_t(\dot B^{\frac  n2+1}_{2,1})}+ \|(\q{v})^h\|_{L^1_t(\dot B^{\frac  n2-1}_{2,1})}\nonumber\\
&\quad\lesssim \|(\ga_1)_0^h\|_{\dot B^{\frac  n2-1}_{2,1}}+\|(\aa_0,\tau_0)^h\|_{\dot B^{\frac  n2}_{2,1}}\nonumber\\
&\quad\quad+\int^t_0(\|(\aa,\tau)^\ell\|_{\dot{B}_{2,1}^{\frac n2-1}}+\|(\aa,\tau)^h\|_{\dot{B}_{2,1}^{\frac n2}})(\| {v}^\ell\|_{\dot B^{\frac  n2}_{2,1}}+\| {v}^h\|_{\dot B^{\frac  n2-1}_{2,1}})\,ds\nonumber\\
&\quad\quad+\int^t_0(\|\aa^\ell\|_{\dot{B}_{2,1}^{\frac n2-1}}+\|\aa^h\|_{\dot{B}_{2,1}^{\frac n2}}+1)(\|\aa^\ell\|_{\dot{B}_{2,1}^{\frac n2-1}}+\|\aa^h\|_{\dot{B}_{2,1}^{\frac n2}})(\|u\|_{\dot{B}_{2,1}^{\frac n2+1}}+\|{v}^\ell\|_{\dot{B}_{2,1}^{\frac n2}})\,ds\nonumber\\
&\quad\quad+ \int^t_0(\big\|(\aa,u,\tau)^\ell\big\|_{\dot{B}^{ \frac n2-1}_{2,1}}+\big\|u^h\big\|_{\dot{B}^{\frac n2-1}_{2,1}}+\big\|(\aa,\tau)^h\big\|_{\dot{B}^{\frac n2}_{2,1}})(\|u^\ell\|_{\dot{B}_{2,1}^{\frac n2+1}}+\|u^h\|_{\dot{B}_{2,1}^{\frac n2+1}})\,ds.
\end{align}
The above estimate is still valid for $\ga_2$ and $\p{v}$:
\begin{align}\label{shiping14}
&\|\ga_2^h\|_{ \widetilde{L}_t^\infty(\dot B^{\frac  n2-1}_{2,1})}+\|(\p{v})^h\|_{ \widetilde{L}_t^\infty(\dot B^{\frac  n2-1}_{2,1})}+ \|\ga_2^h\|_{L^1_t(\dot B^{\frac  n2+1}_{2,1})}+ \|(\p{v})^h\|_{L^1_t(\dot B^{\frac  n2-1}_{2,1})}\nonumber\\
&\quad\lesssim \|(\ga_2)_0^h\|_{\dot B^{\frac  n2-1}_{2,1}}+\|(\aa_0,\tau_0)^h\|_{\dot B^{\frac  n2}_{2,1}}
\nonumber\\
&\quad\quad+\int^t_0(\|(\aa,\tau)^\ell\|_{\dot{B}_{2,1}^{\frac n2-1}}+\|(\aa,\tau)^h\|_{\dot{B}_{2,1}^{\frac n2}})(\| {v}^\ell\|_{\dot B^{\frac  n2}_{2,1}}+\| {v}^h\|_{\dot B^{\frac  n2-1}_{2,1}})\,ds\nonumber\\
&\quad\quad+\int^t_0(\|\aa^\ell\|_{\dot{B}_{2,1}^{\frac n2-1}}+\|\aa^h\|_{\dot{B}_{2,1}^{\frac n2}}+1)(\|\aa^\ell\|_{\dot{B}_{2,1}^{\frac n2-1}}+\|\aa^h\|_{\dot{B}_{2,1}^{\frac n2}})(\|u\|_{\dot{B}_{2,1}^{\frac n2+1}}+\|{v}^\ell\|_{\dot{B}_{2,1}^{\frac n2}})\,ds\nonumber\\
&\quad\quad+ \int^t_0(\big\|(\aa,u,\tau)^\ell\big\|_{\dot{B}^{ \frac n2-1}_{2,1}}+\big\|u^h\big\|_{\dot{B}^{\frac n2-1}_{2,1}}+\big\|(\aa,\tau)^h\big\|_{\dot{B}^{\frac n2}_{2,1}})(\|u^\ell\|_{\dot{B}_{2,1}^{\frac n2+1}}+\|u^h\|_{\dot{B}_{2,1}^{\frac n2+1}})\,ds.
\end{align}
By an elementary manipulation, we have 
\begin{align}
\|u^h\|_{ \widetilde{L}_t^\infty(\dot B^{\frac  n2-1}_{2,1})}
\lesssim&
\|(\ga_1)^h+\Delta^{-1}(\p {v})^h+(\ga_2)^h+\Delta^{-1}(\q {v})^h\|_{ \widetilde{L}_t^\infty(\dot B^{\frac  n2-1}_{2,1})}\nonumber\\
\lesssim&
\|(\ga_1,\ga_2)^h\|_{ \widetilde{L}_t^\infty(\dot B^{\frac  n2-1}_{2,1})}+\|{v}^h\|_{ \widetilde{L}_t^\infty(\dot B^{\frac  n2-1}_{2,1})},\label{shiping16}
\end{align}
\begin{align}
\|u^h\|_{ {L}_t^1(\dot B^{\frac  n2+1}_{2,1})}
\lesssim&
\|(\ga_1,\ga_2)^h\|_{ {L}_t^1(\dot B^{\frac  n2-1}_{2,1})}+\|{v}^h\|_{ {L}_t^1(\dot B^{\frac  n2-1}_{2,1})}.\label{shiping17}
\end{align}

The combination of  \eqref{shiping13}-\eqref{shiping16} with \eqref{shiping17} implies that
\begin{align}\label{shiping18}
&\|u^h\|_{ \widetilde{L}_t^\infty(\dot B^{\frac  n2-1}_{2,1})}+\|{v}^h\|_{ \widetilde{L}_t^\infty(\dot B^{\frac  n2-1}_{2,1})}+ \|u^h\|_{L^1_t(\dot B^{\frac  n2+1}_{2,1})}+ \|{v}^h\|_{L^1_t(\dot B^{\frac  n2-1}_{2,1})}\nonumber\\
&\quad\lesssim \|u_0^h\|_{\dot B^{\frac  n2-1}_{2,1}}+\|(\aa_0,\tau_0)^h\|_{\dot B^{\frac  n2}_{2,1}}
\nonumber\\
&\quad\quad+\int^t_0(\|(\aa,\tau)^\ell\|_{\dot{B}_{2,1}^{\frac n2-1}}+\|(\aa,\tau)^h\|_{\dot{B}_{2,1}^{\frac n2}})(\| {v}^\ell\|_{\dot B^{\frac  n2}_{2,1}}+\| {v}^h\|_{\dot B^{\frac  n2-1}_{2,1}})\,ds\nonumber\\
&\quad\quad+\int^t_0(\|\aa^\ell\|_{\dot{B}_{2,1}^{\frac n2-1}}+\|\aa^h\|_{\dot{B}_{2,1}^{\frac n2}}+1)(\|\aa^\ell\|_{\dot{B}_{2,1}^{\frac n2-1}}+\|\aa^h\|_{\dot{B}_{2,1}^{\frac n2}})(\|u\|_{\dot{B}_{2,1}^{\frac n2+1}}+\|{v}^\ell\|_{\dot{B}_{2,1}^{\frac n2}})\,ds\nonumber\\
&\quad\quad+ \int^t_0(\big\|(\aa,u,\tau)^\ell\big\|_{\dot{B}^{ \frac n2-1}_{2,1}}+\big\|u^h\big\|_{\dot{B}^{\frac n2-1}_{2,1}}+\big\|(\aa,\tau)^h\big\|_{\dot{B}^{\frac n2}_{2,1}})(\|u^\ell\|_{\dot{B}_{2,1}^{\frac n2+1}}+\|u^h\|_{\dot{B}_{2,1}^{\frac n2+1}})\,ds.
\end{align}

Now we show the bound of  $\|(\aa,\tau)^h\|_{ \widetilde{L}_t^\infty(\dot B^{\frac  n2}_{2,1})}$ by using the estimate of $\|u^h\|_{L^1_t(\dot B^{\frac  n2+1}_{2,1})}$ and \eqref{m}.
Indeed,
from  \eqref{m}, we  infer that
\begin{align}\label{shiping20}
&\|\aa^h\|_{ \widetilde{L}_t^\infty(\dot B^{\frac  n2}_{2,1})}
\lesssim \|\aa_0^h\|_{\dot B^{\frac  n2}_{2,1}}+\int_0^t(\|\div u^h\|_{\dot B^{\frac  n2}_{2,1}}+\|(k(a) \div u)^h\|_{\dot B^{\frac  n2}_{2,1}})\,ds
\nonumber\\
&\quad\quad+\int_0^t\|\div u\|_{L^\infty}\|\aa^h\|_{\dot B^{\frac  n2}_{2,1}}ds+\int_0^t\sum_{j\ge j_0}2^{\frac {n}{2}j}\|[\ddj,u\cdot\nabla]\aa\|_{L^2}\,ds\\
&\quad\lesssim \|\aa_0^h\|_{\dot B^{\frac  n2}_{2,1}}+\int_0^t\|u^h\|_{\dot B^{\frac  n2+1}_{2,1}}\,ds
+\int_0^t(\|u^\ell\|_{\dot{B}_{2,1}^{\frac n2+1}}+\|u^h\|_{\dot{B}_{2,1}^{\frac n2+1}})(\|\aa^\ell\|_{\dot{B}_{2,1}^{\frac n2-1}}+\|\aa^h\|_{\dot{B}_{2,1}^{\frac n2}})\,ds\nonumber
\end{align}
and
\begin{align}\label{shiping21}
\|\tau^h\|_{ \widetilde{L}_t^\infty(\dot B^{\frac  n2}_{2,1})}
\lesssim& \|\tau_0^h\|_{\dot B^{\frac  n2}_{2,1}}\!\!+\!\!\int_0^t\|u^h\|_{\dot B^{\frac  n2+1}_{2,1}}\,ds
\nonumber\\
&+\int_0^t(\|u^\ell\|_{\dot{B}_{2,1}^{\frac n2+1}}+\|u^h\|_{\dot{B}_{2,1}^{\frac n2+1}})(\|\tau^\ell\|_{\dot{B}_{2,1}^{\frac n2-1}}+\|\tau^h\|_{\dot{B}_{2,1}^{\frac n2}})\,ds.
\end{align}
Multiplying \eqref{yiming18} by a suitable large constant   and then inserting into  \eqref{shiping20} and \eqref{shiping21}, we get
\begin{align}\label{shiping22}
&\|(u^h,{v}^h)\|_{ \widetilde{L}_t^\infty(\dot B^{\frac  n2-1}_{2,1})}+\|(\aa^h,\tau^h)\|_{ \widetilde{L}_t^\infty(\dot B^{\frac  n2}_{2,1})}+ \|u^h\|_{L^1_t(\dot B^{\frac  n2+1}_{2,1})}+ \|{v}^h\|_{L^1_t(\dot B^{\frac  n2-1}_{2,1})}\nonumber\\
&\quad\lesssim \|u_0^h\|_{\dot B^{\frac  n2-1}_{2,1}}+\|(\aa^h_0,\tau_0^h)\|_{\dot B^{\frac  n2}_{2,1}}
\nonumber\\
&\quad\quad+\int^t_0(\|(\aa,\tau)^\ell\|_{\dot{B}_{2,1}^{\frac n2-1}}+\|(\aa,\tau)^h\|_{\dot{B}_{2,1}^{\frac n2}})(\| {v}^\ell\|_{\dot B^{\frac  n2}_{2,1}}+\| {v}^h\|_{\dot B^{\frac  n2-1}_{2,1}})\,ds\nonumber\\
&\quad\quad+\int^t_0(\|\aa^\ell\|_{\dot{B}_{2,1}^{\frac n2-1}}+\|\aa^h\|_{\dot{B}_{2,1}^{\frac n2}}+1)(\|\aa^\ell\|_{\dot{B}_{2,1}^{\frac n2-1}}+\|\aa^h\|_{\dot{B}_{2,1}^{\frac n2}})(\|u\|_{\dot{B}_{2,1}^{\frac n2+1}}+\|{v}^\ell\|_{\dot{B}_{2,1}^{\frac n2}})\,ds\nonumber\\
&\quad\quad+ \int^t_0(\big\|(\aa,u,\tau)^\ell\big\|_{\dot{B}^{ \frac n2-1}_{2,1}}+\big\|u^h\big\|_{\dot{B}^{\frac n2-1}_{2,1}}+\big\|(\aa,\tau)^h\big\|_{\dot{B}^{\frac n2}_{2,1}})(\|u^\ell\|_{\dot{B}_{2,1}^{\frac n2+1}}+\|u^h\|_{\dot{B}_{2,1}^{\frac n2+1}})\,ds.
\end{align}

Combining \eqref{yiming11} with  \eqref{shiping22}, we finally obtain the  following desired  estimate
\begin{align}\label{shiping24}
&\|(\aa,u,\tau)^\ell\|_{\widetilde{L}_t^{\infty}(\dot{B}_{2,1}^{\frac n2-1})}+\|u^h\|_{ \widetilde{L}_t^\infty(\dot B^{\frac  n2-1}_{2,1})}+\|(\aa^h,\tau^h)\|_{ \widetilde{L}_t^\infty(\dot B^{\frac  n2}_{2,1})}\nonumber\\
&\quad\quad+\|(u,\La^{-1}{v})^\ell\|_{L^1_t(\dot{B}_{2,1}^{\frac n2+1})}+ \|u^h\|_{L^1_t(\dot B^{\frac  n2+1}_{2,1})}+ \|{v}^h\|_{L^1_t(\dot B^{\frac  n2-1}_{2,1})}\nonumber\\
&\quad\lesssim \|( \aa_0,u_0,\tau_0)^\ell\|_{\dot{B}_{2,1}^{\frac {n}{2}-1}}+\|u_0^h\|_{\dot B^{\frac  n2-1}_{2,1}}+\|(\aa^h_0,\tau_0^h)\|_{\dot B^{\frac  n2}_{2,1}}
\nonumber\\
&\quad\quad+\int^t_0(\|(\aa,\tau)^\ell\|_{\dot{B}_{2,1}^{\frac n2-1}}+\|(\aa,\tau)^h\|_{\dot{B}_{2,1}^{\frac n2}}+\|u\|_{\dot{B}_{2,1}^{\frac n2-1}})(\| {v}^\ell\|_{\dot B^{\frac  n2}_{2,1}}+\| {v}^h\|_{\dot B^{\frac  n2-1}_{2,1}})\,ds\nonumber\\
&\quad\quad+\int^t_0(\|\aa^\ell\|_{\dot{B}_{2,1}^{\frac n2-1}}+\|\aa^h\|_{\dot{B}_{2,1}^{\frac n2}}+1)(\|\aa^\ell\|_{\dot{B}_{2,1}^{\frac n2-1}}+\|\aa^h\|_{\dot{B}_{2,1}^{\frac n2}})(\|u\|_{\dot{B}_{2,1}^{\frac n2+1}}+\|{v}^\ell\|_{\dot{B}_{2,1}^{\frac n2}})\,ds\nonumber\\
&\quad\quad+ \int^t_0(\big\|(\aa,u,\tau)^\ell\big\|_{\dot{B}^{ \frac n2-1}_{2,1}}+\big\|u^h\big\|_{\dot{B}^{\frac n2-1}_{2,1}}+\big\|(\aa,\tau)^h\big\|_{\dot{B}^{\frac n2}_{2,1}})(\|u^\ell\|_{\dot{B}_{2,1}^{\frac n2+1}}+\|u^h\|_{\dot{B}_{2,1}^{\frac n2+1}})\,ds.
\end{align}

\subsection{ Completing  the proof of  Theorem \ref{dingli}}
As stated at the beginning of this section, the proof of the  existence of the local solution to \eqref{m} is standard and is thus  omitted. The interested readers may  refer to \cite{chenqionglei}, \cite{fang2014} and  \cite{qian}  for the standard procedure to the local well-posedness.
Therefore,    there existence a constant $T>0$ such that the system \eqref{m} has the required unique and  local solution $(a,u,\tau)$ existing on the time interval  $[0,T)$.
Then, the proof of Theorem \ref{dingli} is reduced to show the extension of $T$ to the infinity,  under the assumption of \eqref{smallness}.

Let
\begin{align*}
X_0\stackrel{\mathrm{def}}{=}&\|( \aa_0,u_0,\tau_0)^\ell\|_{\dot{B}_{2,1}^{\frac {n}{2}-1}}+\|u_0^h\|_{\dot B^{\frac  n2-1}_{2,1}}+\|(\aa^h_0,\tau_0^h)\|_{\dot B^{\frac  n2}_{2,1}},
\end{align*}
and
\begin{align*}
X(t)\stackrel{\mathrm{def}}{=}&\|(\aa,u,\tau)^\ell\|_{\widetilde{L}_t^{\infty}(\dot{B}_{2,1}^{\frac n2-1})}+\|u^h\|_{ \widetilde{L}_t^\infty(\dot B^{\frac  n2-1}_{2,1})}+\|(\aa^h,\tau^h)\|_{ \widetilde{L}_t^\infty(\dot B^{\frac  n2}_{2,1})}\nonumber\\
&\quad\quad+\|(u,\La^{-1}{v})^\ell\|_{L^1_t(\dot{B}_{2,1}^{\frac n2+1})}+ \|u^h\|_{L^1_t(\dot B^{\frac  n2+1}_{2,1})}+ \|{v}^h\|_{L^1_t(\dot B^{\frac  n2-1}_{2,1})}.
\end{align*}
In view of  \eqref{shiping24}, we  get that
\begin{align}\label{xintang}
X(t)\le X_0+C(X(t))^2(1+CX(t)) \,\,\mbox{ for }\,\, t \in [0,T).
\end{align}
By  the  smallness initial data assumption \eqref{smallness},  there exists a positive constant $C_0$ such that
$X_0\leq C_0 \varepsilon$. Hence, 
 for suitable time $T$, we have
\begin{equation}\label{re}
 X (t) \leq 2 C_0\ \varepsilon , \quad   \; t \in [0, T).
\end{equation}
Let $T^{*}$ be the largest possible time of $T$ for the validity of  \eqref{re}.  This is confirmed by \eqref{xintang} and
 the  standard continuation argument  due to the smallness assumption of $\varepsilon$.

Moreover, from the  above argument, we have
\begin{align*}
\|a^\ell\|_{\dot{B}_{2,1}^{\frac n2-1}}\lesssim\|(\sa(\aa))^\ell\|_{\dot{B}_{2,1}^{\frac n2-1}}\lesssim\|\aa^\ell\|_{\dot{B}_{2,1}^{\frac n2-1}}, \quad
\|a^h\|_{\dot{B}_{2,1}^{\frac n2}}\lesssim\|(\sa(\aa))^h\|_{\dot{B}_{2,1}^{\frac n2}}\lesssim\|\aa^h\|_{\dot{B}_{2,1}^{\frac n2}}.
\end{align*}

As
${v}=\nabla \aa -\div \tau,$
we can get
$\p {v}=-\p\div \tau$.
We thus get the smoothing effect, of the incompressible part of $ \div \tau$, expressed as
\begin{align*}
(\Lambda^{-1}\p\div\tau)^\ell\in L^{1}
([0,T];{\dot{B}}_{2,1}^{\frac n2+1}(\R^n)),
 \quad (\Lambda^{-1}\p\div\tau)^h\in L^{1}
([0,T];{\dot{B}}_{2,1}^{\frac n2}(\R^n)),
\end{align*}
for any $T>0$. Consequently, we complete the proof of Theorem \ref{dingli}.

\bigskip
\setcounter{theorem}{0}
\setcounter{equation}{0}

\section{The proof of Theorem \ref{dingli2}}\label{Section4}
In this section, we shall follow the method used in \cite{guoyan} and \cite{xujiang2019arxiv} to get the decay rate of the solution derived  in the previous section.
From the proof of Theorem \ref{dingli} or the derivation of \eqref{shiping24}, we get
 \begin{align}\label{jian18-16}
&\frac{d}{dt}(\|( u,\La^{-1}{v})^\ell\|_{\dot{B}_{2,1}^{\frac n2-1}}+\|u^h\|_{\dot B^{\frac  n2-1}_{2,1}}+\|(\La^{-1}{v})^h\|_{\dot B^{\frac  n2}_{2,1}})\nonumber\\
&\quad\quad\quad+C(\|(u,\La^{-1}{v})^\ell\|_{\dot{B}_{2,1}^{\frac n2+1}}+\|u^h\|_{\dot B^{\frac  n2+1}_{2,1}}+\|(\La^{-1}{v})^h\|_{\dot B^{\frac  n2}_{2,1}})\nonumber\\
&\quad\quad\le C(\|(u,\La^{-1}{v})^\ell\|_{\dot{B}_{2,1}^{\frac n2+1}}+\|u^h\|_{\dot B^{\frac  n2+1}_{2,1}}+\|(\La^{-1}{v})^h\|_{\dot B^{\frac  n2}_{2,1}})\nonumber\\
&\quad\quad\quad\quad\times(\big\|(\aa,u,\tau)^\ell\big\|_{\dot{B}^{ \frac n2-1}_{2,1}}+\big\|\aa^\ell\big\|^2_{\dot{B}^{ \frac n2-1}_{2,1}}+\big\|u^h\big\|_{\dot{B}^{\frac n2-1}_{2,1}}+\big\|(\aa,\tau)^h\big\|_{\dot{B}^{\frac n2}_{2,1}}+\big\|\aa^h\big\|^2_{\dot{B}^{\frac n2}_{2,1}}).
\end{align}
 It follows from  \eqref{xiaonorm} in Theorem \ref{dingli} that
$$\|(\aa,u,\tau)^\ell\|_{\widetilde{L}_t^{\infty}(\dot{B}_{2,1}^{\frac n2-1})}+\|u^h\|_{ \widetilde{L}_t^\infty(\dot B^{\frac  n2-1}_{2,1})}+\|(\aa^h,\tau^h)\|_{ \widetilde{L}_t^\infty(\dot B^{\frac  n2}_{2,1})}
\leq Cc_0\ll 1,$$
for all $t\geq0$.

 Thus, absorbing all the terms on the right-hand  side of  \eqref{jian18-16} to the left, we have
\begin{align}\label{jian18}
&\frac{d}{dt}(\|( u,\La^{-1}{v})^\ell\|_{\dot{B}_{2,1}^{\frac n2-1}}+\|u^h\|_{\dot B^{\frac  n2-1}_{2,1}}+\|(\La^{-1}{v})^h\|_{\dot B^{\frac  n2}_{2,1}})\nonumber\\
&\quad+\bar{c}(\|(u,\La^{-1}{v})^\ell\|_{\dot{B}_{2,1}^{\frac n2+1}}+\|u^h\|_{\dot B^{\frac  n2+1}_{2,1}}+\|(\La^{-1}{v})^h\|_{\dot B^{\frac  n2}_{2,1}})\le 0.
\end{align}

In order to derive the decay estimate of the solution given in Theorem \ref{dingli}, we need  to get  a Lyapunov-type differential inequality from \eqref{jian18}. This inequality can be obtained from interpolation inequality, which heavily relies on the  bound 
\begin{equation}\label{new3}
\|(u,\La^{-1}{v})\|^{\ell}_{\dot{B}^{\sigma}_{2,1}}\le C, \,\,\hbox{  $-\frac n2<\sigma<\frac n2-1$}.
\end{equation}
Since ${v}=\nabla \aa -\div \tau$, we only need to estimate  $\|(\aa,u,\tau)\|^{\ell}_{\dot{B}^{\sigma}_{2,1}}$.
Due to the linear function assumption of $P(\rho)$, we have 
$\aa=\tilde{c} a$ for a positive constant $\tilde{c}.$ Hence,  the following argument aims to derive  the bound
\begin{equation}\|(a,u,\tau)\|^{\ell}_{\dot{B}^{\sigma}_{2,1}}\le C, \,\,\hbox{  $-\frac n2<\sigma<\frac n2-1$}.
\end{equation}

From \eqref{m} and the  derivation of \eqref{ert0}  it follows that
\begin{align}\label{jian1}
\|(a,u,\tau)^\ell\|_{\widetilde{L}_t^{\infty}(\dot{B}_{2,1}^\sigma)}
\lesssim&\|(a_0,u_0,\tau_0)^\ell\|_{\dot{B}_{2,1}^{\frac {n}{2}-1}}
+ \int^t_0\|\nabla  u\|_{L^\infty}\|(a,u,\tau)^\ell\|_{\dot{B}_{2,1}^\sigma}\,ds\nonumber\\
&+ \int^t_0\|(a \div u)^\ell\|_{\dot{B}_{2,1}^\sigma}\,ds+ \int^t_0\|(g(\tau, \nabla u))^\ell\|_{\dot{B}_{2,1}^\sigma}\,ds
\nonumber\\
&+ \int^t_0\|(I(a){v})^\ell\|_{\dot{B}_{2,1}^\sigma}\,ds+ \int^t_0\|(I(a)(\Delta u+\nabla\div u))^\ell\|_{\dot{B}_{2,1}^\sigma}\,ds\nonumber\\
&+\int^t_0\sum_{j\le j_0}2^{\sigma j}\big(\| [\ddj ,u\cdot \na ] ( a,u,\tau)\|_{L^2}\big)\,ds.
\end{align}
By virtue of embedding relation $\dot{B}_{2,1}^{\frac n2}(\R^n)\hookrightarrow L^\infty(\R^n),$ we have
\begin{align}\label{jian2}
\|\nabla  u\|_{L^\infty}\|(a,u,\tau)^\ell\|_{\dot{B}_{2,1}^\sigma}
\lesssim&\|\nabla  u\|_{\dot{B}_{2,1}^{\frac n2}}\|(a,u,\tau)^\ell\|_{\dot{B}_{2,1}^\sigma}\nonumber\\
\lesssim&(\|u^\ell\|_{\dot{B}_{2,1}^{\frac n2+1}}+\|  u^h\|_{\dot{B}_{2,1}^{\frac n2+1}})\|(a,u,\tau)^\ell\|_{\dot{B}_{2,1}^\sigma}.
\end{align}

In the further argument,  we shall use repeatedly  the following  product law:
\begin{align}\label{jian5}
\|(bc^h)^\ell\|_{\dot{B}_{2,1}^{\sigma}}\lesssim\|b\|_{\dot{B}_{2,1}^{\frac n2-1}}\|c^h\|_{\dot{B}_{2,1}^{\frac n2-1}},\quad\hbox{for any $-\frac n2<\sigma<\frac n2-1$}.
\end{align}

Indeed, by
Bony's decomposition:
$$
bc^h= \dot{T}_{c^h}b+\dot{R}(c^h,b)+\dot{T}_{b}c^h.
$$
By  Lemma \ref{fangji}, the first term on the right-hand side of the previous equation  is  estimated as
\begin{align}\label{jiang998}
\|\dot{T}_{c^h}b\|_{\dot B^{\sigma}_{2,1}}^\ell
\lesssim&\|\dot{T}_{c^h}b\|_{\dot B^{-\frac n2}_{2,1}}^\ell
\lesssim \|c^h\|_{\dot B^{1-n}_{\infty,\infty}}\|b\|_{\dot B^{\frac n2-1}_{2,1}}\nonumber\\
\lesssim& \|c^h\|_{\dot B^{1-\frac n2}_{2,1}}\|b\|_{\dot B^{\frac n2-1}_{2,1}}
\lesssim \|c^h\|_{\dot B^{\frac n2-1}_{2,1}}\|b\|_{\dot B^{\frac n2-1}_{2,1}},
\end{align}
where we have used the high  frequency property of $c^h$ and the  fact $1-\frac n2\le\frac n2-1$ in the last inequality of the previous equation.

 Similarly, by the low  frequency property and the condition $-\frac n2< \sigma<\frac n2-1$, the remaining terms are bounded as
 \begin{align}\label{jiang999}
\|\dot{T}_b{c^h}\|_{\dot B^{\sigma}_{2,1}}^\ell
\lesssim&\|\dot{T}_b{c^h}\|_{\dot B^{-\frac n2}_{2,1}}^\ell
\lesssim \|b\|_{\dot B^{-1}_{\infty,\infty}}\|c^h\|_{\dot B^{1-\frac n2}_{2,1}}\nonumber\\
\lesssim& \|b\|_{\dot B^{\frac n2-1}_{2,1}}\|c^h\|_{\dot B^{1-\frac n2}_{2,1}}
\lesssim \|b\|_{\dot B^{\frac n2-1}_{2,1}}\|c^h\|_{\dot B^{\frac n2-1}_{2,1}}.
\end{align}
and
\begin{align}\label{jiang99899}
\|\dot{R}(c^h,b)\|^\ell_{\dot B^{\sigma}_{2,1}}
\lesssim&\|\dot{R}(c^h,b)\|_{\dot B^{-\frac n2}_{2,\infty}}\lesssim\|\dot{R}(c^h,b)\|_{\dot B^{0}_{1,\infty}}\nonumber\\
\lesssim& \|c^h\|_{\dot B^{1-\frac n2}_{2,1}}\|b\|_{\dot B^{\frac n2-1}_{2,1}}
\lesssim \|c^h\|_{\dot B^{\frac n2-1}_{2,1}}\|b\|_{\dot B^{\frac n2-1}_{2,1}},
\end{align}
where we have used  Lemma \ref{fangji}.
Hence \eqref{jian5} follows from \eqref{jiang998}-\eqref{jiang99899}.

Thanks to Lemma \ref{daishu}, we  infer  that
\begin{align}\label{jian3}
\|(a \div u)^\ell\|_{\dot{B}_{2,1}^\sigma}
\lesssim& \|a\|_{\dot{B}_{2,1}^\sigma}\|\div u\|_{\dot{B}_{2,1}^{\frac n2}}\nonumber\\
\lesssim& (\|a^\ell\|_{\dot{B}_{2,1}^\sigma}+\|a^h\|_{\dot{B}_{2,1}^{\sigma}})\| u\|_{\dot{B}_{2,1}^{\frac n2+1}}\nonumber\\
\lesssim& (\|a^\ell\|_{\dot{B}_{2,1}^\sigma}+\|a^h\|_{\dot{B}_{2,1}^{\frac n2}})(\|u^\ell\|_{\dot{B}_{2,1}^{\frac n2+1}}+\|  u^h\|_{\dot{B}_{2,1}^{\frac n2+1}}).
\end{align}

In  what follows, we  shall   repeatedly use  the following fact:
\begin{align}\label{airen}
\|a^h\|_{\dot{B}_{2,1}^{\sigma}}\lesssim\|a^h\|_{\dot{B}_{2,1}^{\frac n2}},\quad\hbox{ $-\frac n2<\sigma<\frac n2-1$}.
\end{align}

Following the  derivation of \eqref{jian3}, we have
\begin{align}\label{jian4}
\| (g(\tau, \nabla u))^\ell\|_{\dot{B}_{2,1}^\sigma}
\lesssim&(\|\tau^\ell\|_{\dot{B}_{2,1}^\sigma}+\|\tau^h\|_{\dot{B}_{2,1}^{\frac n2}})(\|u^\ell\|_{\dot{B}_{2,1}^{\frac n2+1}}+\|  u^h\|_{\dot{B}_{2,1}^{\frac n2+1}}).
\end{align}
For  the estimate of the nonlinear term $\|(I(a) {v})^\ell\|_{\dot{B}_{2,1}^\sigma}$, 
we use  Lemma \ref{daishu} and \eqref{jian5} to produce that
\begin{align}\label{jian6}
\|(I(a) {v})^\ell\|_{\dot{B}_{2,1}^\sigma}
\lesssim& \|I(a){v}^\ell\|_{\dot{B}_{2,1}^\sigma}+\|I(a){v}^h\|_{\dot{B}_{2,1}^{\sigma}}\nonumber\\
\lesssim& \|I(a)\|_{\dot{B}_{2,1}^\sigma}\|{v}^\ell\|_{\dot{B}_{2,1}^{\frac n2}}+\|I(a)\|_{\dot{B}_{2,1}^{\frac n2-1}}\|{v}^h\|_{\dot{B}_{2,1}^{\frac n2-1}}.
\end{align}
Recall that $I(a)=\frac{a}{1+a}=a-aI(a)$. Hence, we have
\begin{align}\label{jian8}
\|I(a)\|_{\dot{B}_{2,1}^\sigma}\lesssim& \|a-aI(a)\|_{\dot{B}_{2,1}^\sigma}
\lesssim (1+\|I(a)\|_{\dot{B}_{2,1}^{\frac n2}})\|a\|_{\dot{B}_{2,1}^\sigma}\nonumber\\
\lesssim& (1+\|a\|_{\dot{B}_{2,1}^{\frac n2}})(\|a^\ell\|_{\dot{B}_{2,1}^\sigma}+\|a^h\|_{\dot{B}_{2,1}^{\frac n2}})\nonumber\\
\lesssim&(1+\|a^\ell\|_{\dot{B}_{2,1}^{\frac n2-1}}+\|a^h\|_{\dot{B}_{2,1}^{\frac n2}})\|a^\ell\|_{\dot{B}_{2,1}^\sigma}
+(1+\|a^\ell\|_{\dot{B}_{2,1}^{\frac n2-1}}+\|a^h\|_{\dot{B}_{2,1}^{\frac n2}})
\|a^h\|_{\dot{B}_{2,1}^{\frac n2}}
\end{align}
and, similarly,
\begin{align}\label{jian9}
\|I(a)\|_{\dot{B}_{2,1}^{\frac n2-1}}
\lesssim&\|a-aI(a)\|_{\dot{B}_{2,1}^{\frac n2-1}}\nonumber\\
\lesssim&\|a\|_{\dot{B}_{2,1}^{\frac n2-1}}+\|a\|_{\dot{B}_{2,1}^{\frac n2-1}}\|I(a)\|_{\dot{B}_{2,1}^{\frac n2}}\nonumber\\
\lesssim&\|a\|_{\dot{B}_{2,1}^{\frac n2-1}}+\|a\|_{\dot{B}_{2,1}^{\frac n2-1}}\|a\|_{\dot{B}_{2,1}^{\frac n2}}\nonumber\\
\lesssim&(\|a^\ell\|_{\dot{B}_{2,1}^{\frac n2-1}}+\|a^h\|_{\dot{B}_{2,1}^{\frac n2}})(1+\|a^\ell\|_{\dot{B}_{2,1}^{\frac n2-1}}+\|a^h\|_{\dot{B}_{2,1}^{\frac n2}}).
\end{align}
Inserting  \eqref{jian8} and \eqref{jian9} into \eqref{jian6}, we have
\begin{align}\label{jian10}
\|(I(a) {v})^\ell\|_{\dot{B}_{2,1}^\sigma}
\lesssim&(1+\|a^\ell\|_{\dot{B}_{2,1}^{\frac n2-1}}+\|a^h\|_{\dot{B}_{2,1}^{\frac n2}})\|{v}^\ell\|_{\dot{B}_{2,1}^{\frac n2}}\|a^\ell\|_{\dot{B}_{2,1}^\sigma}\nonumber\\
&+(1+\|a^\ell\|_{\dot{B}_{2,1}^{\frac n2-1}}+\|a^h\|_{\dot{B}_{2,1}^{\frac n2}})\|a^h\|_{\dot{B}_{2,1}^{\frac n2}}\|{v}^\ell\|_{\dot{B}_{2,1}^{\frac n2}}
\nonumber\\
&+(\|a^\ell\|_{\dot{B}_{2,1}^{\frac n2-1}}+\|a^h\|_{\dot{B}_{2,1}^{\frac n2}})(1+\|a^\ell\|_{\dot{B}_{2,1}^{\frac n2-1}}+\|a^h\|_{\dot{B}_{2,1}^{\frac n2}})\|{v}^h\|_{\dot{B}_{2,1}^{\frac n2-1}}.
\end{align}
and, correspondingly, %
\begin{align}\label{jian11}
&\|(I(a)(\Delta u+\nabla\div u))^\ell\|_{\dot{B}_{2,1}^\sigma}\nonumber\\
&\quad\lesssim(\|a^\ell\|_{\dot{B}_{2,1}^{\frac n2-1}}+\|a^h\|_{\dot{B}_{2,1}^{\frac n2}}+1)\|u^\ell\|_{\dot{B}_{2,1}^{\frac n2+1}}\|a^\ell\|_{\dot{B}_{2,1}^\sigma}\nonumber\\
&\quad\quad+(1+\|a^\ell\|_{\dot{B}_{2,1}^{\frac n2-1}}+\|a^h\|_{\dot{B}_{2,1}^{\frac n2}})\|a^h\|_{\dot{B}_{2,1}^{\frac n2}}\|u^\ell\|_{\dot{B}_{2,1}^{\frac n2+1}}
\nonumber\\
&\quad\quad+(\|a^\ell\|_{\dot{B}_{2,1}^{\frac n2-1}}+\|a^h\|_{\dot{B}_{2,1}^{\frac n2}})(1+\|a^\ell\|_{\dot{B}_{2,1}^{\frac n2-1}}+\|a^h\|_{\dot{B}_{2,1}^{\frac n2}})\|u^h\|_{\dot{B}_{2,1}^{\frac n2+1}}.
\end{align}

With the aid of Lemma \ref{jiaohuanzi}, we have
\begin{align}\label{jian12}
\sum_{j\le j_0}2^{\sigma j}\| [\ddj ,u\cdot \na ]   (a,u,\tau)\|_{L^2}
\lesssim&\big\|\nabla u\big\|_{\dot{B}^{\frac{n}{2}}_{2,1}}\big\|(a,u,\tau)\big\|_{\dot{B}^{ \sigma}_{2,1}}\nonumber\\
\lesssim&(\|u^\ell\|_{\dot{B}_{2,1}^{\frac n2+1}}+\|  u^h\|_{\dot{B}_{2,1}^{\frac n2+1}})(\big\|u^h\big\|_{\dot{B}^{\frac n2-1}_{2,1}}+\big\|(a,\tau)^h\big\|_{\dot{B}^{\frac n2}_{2,1}})
\nonumber\\
&+(\|u^\ell\|_{\dot{B}_{2,1}^{\frac n2+1}}+\|  u^h\|_{\dot{B}_{2,1}^{\frac n2+1}})\big\|(a,u,\tau)^\ell\big\|_{\dot{B}^{ \sigma}_{2,1}}.
\end{align}
Inserting
\eqref{jian2}, \eqref{jian3}, \eqref{jian4}, \eqref{jian10}--\eqref{jian12} into \eqref{jian1} gives
\begin{align}\label{jian13}
&\|(a,u,\tau)^\ell\|_{\widetilde{L}_t^{\infty}(\dot{B}_{2,1}^\sigma)}
\lesssim\|(a_0,u_0,\tau_0)^\ell\|_{\dot{B}_{2,1}^{\frac {n}{2}-1}}+\int_0^tG_2(s)ds
+\int_0^tG_1(s)\|(a,u,\tau)^\ell\|_{\dot{B}_{2,1}^\sigma}ds
\end{align}
in which
\begin{align*}
G_1(t)=&(\|a^\ell\|_{\dot{B}_{2,1}^{\frac n2-1}}+\|a^h\|_{\dot{B}_{2,1}^{\frac n2}}+1)(\|u^\ell\|_{\dot{B}_{2,1}^{\frac n2+1}}+\|  u^h\|_{\dot{B}_{2,1}^{\frac n2+1}}+\|{v}^\ell\|_{\dot{B}_{2,1}^{\frac n2}}),\nonumber\\
G_2(t)=&(\|a^\ell\|_{\dot{B}_{2,1}^{\frac n2-1}}+\|a^h\|_{\dot{B}_{2,1}^{\frac n2}}+1)\|(a,\tau)^h\|_{\dot{B}_{2,1}^{\frac n2}}(\|u^\ell\|_{\dot{B}_{2,1}^{\frac n2+1}}+\|{v}^\ell\|_{\dot{B}_{2,1}^{\frac n2}})\nonumber\\
&+(\|a^\ell\|_{\dot{B}_{2,1}^{\frac n2-1}}+\|(a,\tau)^h\|_{\dot{B}_{2,1}^{\frac n2}})(\|a^\ell\|_{\dot{B}_{2,1}^{\frac n2-1}}+\|a^h\|_{\dot{B}_{2,1}^{\frac n2}}+1)(\|u^h\|_{\dot{B}_{2,1}^{\frac n2+1}}+\|{v}^h\|_{\dot{B}_{2,1}^{\frac n2-1}}).
\end{align*}
As $(a,u,\tau)$ is the global solution constructed in Theorem \ref{dingli},  we  get from \eqref{xiaonorm} that
\begin{align}\label{jian15}
\int_0^tG_1(s)ds+\int_0^tG_2(s)ds\le C,
\end{align}
from which and the Gronwall inequality applied to \eqref{jian13}, we have
\begin{align}\label{jian16}
&\|(a,u,\tau)^\ell\|_{\dot{B}^{\sigma}_{2,1}}
\le C,\quad  -\frac n2< \sigma<\frac n2-1,
\end{align}
with $C$  a  constant depending on $n$, $a_0,u_0,$ and $\tau_0$.

Moreover,
from the definition $v=\nabla a-\div\tau$, we  get
\begin{align}\label{jian17}
&\|(\La^{-1}{v})^\ell\|_{\dot{B}^{\sigma}_{2,1}}\le C\|(a,\tau)^\ell\|_{\dot{B}^{\sigma}_{2,1}}
\le C,\quad  -\frac n2< \sigma<\frac n2-1.
\end{align}
Now,
for any $ -\frac n2< \sigma<\frac n2-1$,
it follows from interpolation inequality in Lemma \ref{qianru}  that
\begin{align*}
\|( u,{\La^{-1}{v}})\|^\ell_{\dot{B}_{2,1}^{\frac n2-1}}
\le& C \big(\|( u,{\La^{-1}{v}})\|^\ell_{\dot{B}_{2,1}^{\sigma}}\big)^{\theta_{1}}\big(\|( u,{\La^{-1}{v}})\|^\ell_{\dot{B}_{2,1}^{\frac n2+1}}\big)^{1-\theta_{1}}\nonumber\\
\le&C\big(\|( u,{\La^{-1}{v}})\|^\ell_{\dot{B}_{2,1}^{\frac n2+1}}\big)^{1-\theta_{1}},
\quad \theta_1=\frac{4}{n-2\sigma+2}\in(0,1).
\end{align*}
This  implies that
\begin{align}\label{jian20}
\|( u,{\La^{-1}{v}})\|^\ell_{\dot{B}_{2,1}^{\frac n2+1}}\ge  C\big(\|( u,{\La^{-1}{v}})\|^\ell_{\dot{B}_{2,1}^{\frac n2-1}}\big)^{\frac{1}{1-\theta_{1}}}.
\end{align}
Due to the embedding relations in high frequencies, we deduce from  \eqref{xiaonorm}
that
\begin{align}\label{jian21}
\|u\|^h_{\dot{B}_{2,1}^{\frac n2+1}}\ge C\big(\|u\|^h_{\dot{B}_{2,1}^{\frac n2-1}}\big)^{\frac{1}{1-\theta_{1}}}, \quad\|{\La^{-1}{v}}\|^h_{\dot{B}_{2,1}^{\frac n2}}\ge C\big(\|{\La^{-1}{v}}\|^h_{\dot{B}_{2,1}^{\frac n2}}\big)^{\frac{1}{1-\theta_{1}}}.
\end{align}

Thus, substituting \eqref{jian20} and \eqref{jian21} into \eqref{jian18} yields
\begin{align*}
&\frac{d}{dt}\big(\|(u,{\La^{-1}{v}})\|^\ell_{\dot{B}_{2,1}^{\frac n2-1}}+\|u\|^h_{\dot B^{\frac  n2-1}_{2,1}}+\|{\La^{-1}{v}}\|^h_{\dot B^{\frac  n2}_{2,1}}\big)\\
&\quad
+\bar{c}\big(\|(u,{\La^{-1}{v}})\|^\ell_{\dot{B}_{2,1}^{\frac n2-1}}+\|u\|^h_{\dot B^{\frac  n2-1}_{2,1}}+\|{\La^{-1}{v}}\|^h_{\dot B^{\frac  n2}_{2,1}}\big)^{\frac{ n-2\sigma+2}{ n-2\sigma-2}}\le 0.
\end{align*}
Solving this differential inequality directly, we obtain
\begin{align*}
\|(u,{\La^{-1}{v}})\|^\ell_{\dot{B}_{2,1}^{\frac n2-1}}+\|u\|^h_{\dot B^{\frac  n2-1}_{2,1}}+\|{\La^{-1}{v}}\|^h_{\dot B^{\frac  n2}_{2,1}}\le&C ({\mathcal X}_0^{-\frac{4}{n-2\sigma-2}}+\frac{4\bar{c}}{n-2\sigma-2}t)^{-\frac{n-2\sigma-2}{4}}\\
\le& C(1+t)^{-\frac{n-2\sigma-2}{4}}.
\end{align*}
Moreover,  we further get
\begin{align}\label{jian23}
\|(u,{\La^{-1}{v}})\|_{\dot{B}_{2,1}^{\frac n2-1}}\le C(\|(u,{\La^{-1}{v}})\|^\ell_{\dot{B}_{2,1}^{\frac n2-1}}+\|u\|^h_{\dot{B}_{2,1}^{\frac n2-1}}+\|{\La^{-1}{v}}\|^h_{\dot{B}_{2,1}^{\frac n2}})\le C(1+t)^{-\frac{n-2\sigma-2}{4}}.
\end{align}
For any $\sigma<\gamma<\frac n2-1,$  we get from the embedding relations in the high frequencies that
\begin{align}\label{jian099}
\|(u,{\La^{-1}{v}})^h\|_{\dot{B}_{2,1}^{\gamma}}\le C(\|u\|^h_{\dot{B}_{2,1}^{\frac n2-1}}+\|{\La^{-1}{v}}\|^h_{\dot{B}_{2,1}^{\frac n2}})\le C(1+t)^{-\frac{n-2\sigma-2}{4}},
\end{align}
On the other hand,
by the interpolation inequality, we have
\begin{align*}
\|(u,{\La^{-1}{v}})\|^\ell_{\dot{B}_{2,1}^{\gamma}}
\le&C\big(\|(u,{\La^{-1}{v}})\|^\ell_{\dot{B}_{2,1}^{\sigma}}\big)^{\theta_{2}} \big(\|(u,{\La^{-1}{v}})\|^\ell_{\dot{B}_{2,1}^{\frac n2-1}}\big)^{1-\theta_{2}},\quad \theta_{2}=\frac{\frac n2 -1-\gamma}{\frac n2-1-\sigma}\in (0,1),
\end{align*}
which is combined with  \eqref{jian16}, \eqref{jian17} and  \eqref{jian23} to produce
\begin{align}\label{jian24}
\|(u,{\La^{-1}{v}})\|^\ell_{\dot{B}_{2,1}^{\gamma}}
\le C(1+t)^{-\frac{(\frac n2-\sigma-1)\theta_{2}}{2}}
=C(1+t)^{-\frac{\gamma-\sigma}{2}}.
\end{align}
The combination of \eqref{jian099} and \eqref{jian24} gives
\begin{align*}
\|(u,{\La^{-1}{v}})\|_{\dot{B}_{2,1}^{\gamma}}
\le&C(\|(u,{\La^{-1}{v}})\|^\ell_{\dot{B}_{2,1}^{\gamma}}+\|(u,{\La^{-1}{v}})\|^h_{\dot{B}_{2,1}^{\gamma}})\nonumber\\
\le& C(1+t)^{-\frac{\gamma-\sigma}{2}}.
\end{align*}

Thanks to the embedding relation
$\dot{B}^{0}_{2,1}(\R^n)\hookrightarrow L^2(\R^n)$, we infer that
\begin{align*}
\|\Lambda^{\gamma} (u,{\La^{-1}{v}})\|_{L^2}
\le& C(1+t)^{-\frac{\gamma-\sigma}{2}}.
\end{align*}
For
$2\leq q\leq\infty$, $\frac nq-\frac n2+\sigma<\alpha \leq\frac nq-1$,  $m=\frac n2-1,$  integer $k\ge 0$ and
$$
k\theta_{3}+m(1-\theta_{3})=\alpha+n\Bigl(\frac12-\frac1q\Bigr), \quad
$$
we use  the Gagliardo-Nirenberg type interpolation inequality \cite[Chap. 2]{bcd} to get
\begin{align*}
\|\Lambda^{\alpha}(u,{\La^{-1}{v}})\|_{L^{q}} \le &C\|\Lambda^{m}(u,{\La^{-1}{v}})\|_{L^2}^{1-\theta_{3}}\|\Lambda^{k}(u,{\La^{-1}{v}})\|^{\theta_{3}}_{L^2}
\nonumber\\ \le & C\Big\{(1+t)^{-\frac{m-\sigma}{2}}\Big\}^{1-\theta_{3}}
\Big\{(1+t)^{-\frac{k-\sigma}{2}}\Big\}^{\theta_{3}}
\nonumber\\=& C(1+t)^{-\frac n2(\frac 12-\frac 1q)-\frac {\alpha-\sigma}{2}}.
\end{align*}
Consequently, we  complete the proof of Theorem \ref{dingli2}. \quad\quad$\Box$

\

\noindent {\bf Acknowledgement.}
 This work is supported by NSFC in China  under grant numbers 11601533 and  11571240.

\noindent$\large\mathbf{Conflict\ of\ Interest:}$ The authors declare that they have no conflict of interest.

This research dose not involved human participants and   this article does not
contain any studies with animals performed by any of the authors.

\bigskip
\bigskip


\end{document}